\documentclass[onefignum,onetabnum]{siamonline220329}
\usepackage{algpseudocode}
\usepackage[caption=false]{subfig}
\usepackage{amssymb}
\usepackage{multirow}
\usepackage{booktabs}

\usepackage{array}
\newcolumntype{L}{>{\centering\arraybackslash}m{3cm}}

\algrenewcommand\algorithmicrequire{\textbf{Input:}}
\algrenewcommand\algorithmicensure{\textbf{Output:}}

\DeclareMathOperator*{\argmin}{argmin}

\newcommand{\Vector}{\textbf{vec}}
\newcommand{\Arrayop}{\textbf{array}}

\newcommand{\bfb}{\mathbf b}

\newcommand{\bff}{\mathbf f}

\newcommand{\bfw}{\mathbf w}
\newcommand{\bfz}{\mathbf z}

\newcommand{\bfbtrue}{\bfb^{\mathrm{true}}}

\newcommand{\bfeta}{\boldsymbol{\eta}}

\newcommand{\bfs}{\mathbf s}
\newcommand{\bfq}{\mathbf q}

\newcommand{\bfr}{\mathbf r}

\newcommand{\bfu}{\mathbf u}
\newcommand{\bfv}{\mathbf v}
\newcommand{\bfx}{\mathbf x}
\newcommand{\bfX}{\mathbf X}
\newcommand{\bfy}{\mathbf y}

\newcommand{\bfg}{\mathbf g}

\newcommand{\bfh}{\mathbf h}
\newcommand{\bfd}{\mathbf d}

\newcommand{\bfc}{\mathbf c}

\newcommand{\bfzeta}{\boldsymbol{\zeta}}

\newcommand{\bfxtrue}{\bfx^{\mathrm{true}}}

\newcommand*\samethanks[1][\value{footnote}]{\footnotemark[#1]}

\newtheorem{exmp}{Example}[section]

\makeatother

\title{Split Bregman Isotropic and Anisotropic Image Deblurring with Kronecker Product Sum Approximations using Single Precision Enlarged-GKB or RSVD Algorithms to provide low rank truncated SVDs 
}
   \author{Abdulmajeed Alsubhi{\thanks{School of Mathematical and Statistical Sciences, Arizona State University, Tempe, AZ (\href{mailto:ahalsubh@asu.edu}{ahalsubh@asu.edu})\and (\href{mailto:renaut@asu.edu}{renaut@asu.edu})} \thanks{Department of Mathematics, Faculty of Science, Islamic University of Madinah, Madinah, Saudi Arabia}}\and Rosemary A Renaut{\samethanks[1]}}

\headers{$\ell_1$ Regularized Problems via KP approximation}{Alsubhi and  Renaut}

\begin{document} 

\maketitle

\begin{abstract}
We consider the solution of the $\ell_1$ regularized  image deblurring problem  using isotropic and anisotropic regularization implemented with  the split Bregman algorithm. For  large scale problems, we replace the system matrix $A$ using a Kronecker product approximation obtained via an approximate truncated singular value decomposition for the reordered matrix $\mathcal{R}(A)$. To obtain the approximate decomposition for $\mathcal{R}(A)$ we propose the enlarged Golub Kahan Bidiagonalization algorithm that proceeds by enlarging the Krylov subspace beyond either a given rank for the desired approximation, or uses an automatic stopping test that provides a suitable rank for the approximation. The resultant expansion is contrasted with the use of the truncated and the randomized singular value decompositions with the same number of terms. To further extend the scale of problem that can be considered we implement the determination of the approximation using single precision, while performing all steps for the regularization in standard double precision. The reported numerical tests demonstrate the effectiveness of applying the approximate single precision Kronecker product expansion for $A$, combined with either isotropic or anisotropic regularization implemented using the split Bregman algorithm, for the solution of image deblurring  problems.  As the size of the problem increases, our results demonstrate that  the major costs are associated with determining the Kronecker product approximation, rather than with the cost of the regularization algorithm. Moreover, the enlarged Golub Kahan Bidiagonalization algorithm  competes favorably with the randomized singular value decomposition for estimating the approximate singular value decomposition.
\end{abstract}

\begin{keywords}   $\ell_1$ regularization, Split Bregman, Anisotropic and Isotropic, Kronecker Product, Single and double precision, Golub Kahan Bidiagonalization
\end{keywords}

\begin{AMS} 65F22, 65F10, 68W40
\end{AMS}

\section{Introduction}
We are concerned with the solution of the ill-posed inverse problem
\begin{equation}\label{eq:forward problem}
    A\bfx \approx \bfb = \bfbtrue+ \bfeta,
\end{equation}
where $A\in \mathbb{R}^{M\times N}$,  is a non separable matrix, $\bfx$ is an unknown solution, $\bfb$ is observed data, and $\bfeta$ is a noise vector which is assumed to be Gaussian white noise. The matrix $A$ is ill-conditioned and the solution $\bfx$ of \cref{eq:forward problem} is sensitive to the noise $\bfeta$. Determining the solution $\bfx$ is regarded as an ill-posed problem \cite[section 1.2]{hansen2010discrete}, \cite[section 2.1]{vogel2002computational}.  In particular, the Picard condition, the requirement that the absolute values of the coefficients of the solution go to zero faster than the singular values of $A$, does not hold. Due to the introduction of the noise, the coefficients of $\bfb$ level out at the noise level, while the singular values continue to decay to the machine precision. Consequently, the ratios of the coefficients to the singular values blow up in absolute value and contaminate the desired solution $\bfx$. Without extra consideration, the noise in the measurements will contaminate the solution. Hence, predicting $\bfx$ from $\bfb$ presents algorithmic challenges.

A solution of an ill-posed problem  can be obtained by applying a regularization technique to mitigate the impact of the noise in the measured data \cite{bjorck1996numerical,golub2013matrix}. Specifically, the ill-posed problem is replaced with a stabilizing  well-posed problem. The resulting problem has better numerical properties and, therefore,  may yield a more accurate approximation to the desired solution. A regularized solution can be obtained, for example, by minimizing the Tikhonov problem \cite{fuhry2012new,golub1999tikhonov}. A least squares solution of the Tikhonov regularized problem  tends to be smooth and not able to capture all the details needed for the desired solution. Here we suppose an acceptable solution may be obtained from the solution of the $\ell_1$ regularized inverse problem
\begin{equation}\label{eq:Tv}
    \bfx= \argmin_\bfx\{\|A \bfx - \bfb \|^2_2 + \beta \| L \bfx\|_1\},
\end{equation}
where $\beta$ is a regularization parameter that trades off the contribution of the fidelity term $\|A \bfx - \bfb \|^2_2$ and the regularization term $\| L \bfx\|_1$, and $L \in \mathbb{R}^{P\times N}$ is a regularization matrix. 

Total variation (TV) regularization, as a form of $\ell_1$ regularization,  was first applied to image restoration in \cite{rudin1992nonlinear}. There it was  demonstrated that TV regularization can be used to preserve sharp edges in the desired image. The optimization problem \cref{eq:Tv} depends on $L$ and $\beta$ but no explicit formula for the solution is available, as it is for the Tikhonov problem. 
The appearance of the $\ell_1$ regularization term means that it is more difficult to solve \cref{eq:Tv} than the same problem with the $\ell_1$ term replaced by a $\ell_2$ term, namely the Tikhonov problem. There are many proposed numerical algorithms \cite{chambolle2004algorithm,goldstein2009split,rudin1992nonlinear,vogel1996iterative} that can be used for the solution of \cref{eq:Tv}. Here we focus on an efficient  Split Bregman (SB) algorithm \cite{goldstein2009split} that has been  shown to improve the reconstruction of blurred images, but  requires solving a linear system of equations at each iteration. Thus, achieving an iterative solution with the SB method for large $A$ becomes computationally costly or infeasible unless $A$ has a special structure. For example, when considering an invariant blur for an image with zero boundary conditions, $A$ exhibits a block-Toeplitz-Toeplitz block (BTTB) structure, and the forward and transpose operations can be implemented efficiently using the Fast Fourier transform in two dimensions, e.g.  \cite{hansen2006deblurring,vogel2002computational}.

 When $A$ is large it may be more efficient  to use the Kronecker product (KP) sum approximation that was introduced by Loan and Pitsianis in \cite{loan1992approximation}. This approach requires a reordering of the matrix $A$, yielding a matrix $\mathcal{R}(A)$, for which its SVD leads to the sum approximation of $A$. Computing the SVD of $\mathcal{R}(A)$, however, might also be computationally infeasible because $\mathcal{R}(A)$ inherits the large size of $A$. An alternative approach for constructing a KP decomposition was discussed for banded Toplitz matrices in \cite{nagy1996decomposition}, for structured blur operators in \cite{kamm2000optimal},  and for specific types of boundary conditions in \cite{nagy2003kronecker}.

Under the assumption that $A$ is large, unstructured, and inherently not separable, we focus on the use of the method of Loan and Pitsianis in \cite{loan1992approximation}. This decomposition was extended in  \cite{garvey2018singular} to provide a TSVD approximation to a non-separable operator matrix, but relies on reordering the singular values of the first term in the KP sum to find a truncated expansion to $A$.
Here we propose two alternative approaches based on two algorithms for finding a good low rank TSVD of a matrix:  specifically (i) the well-known randomized singular value decomposition (RSVD) and (ii) the use of  an enlarged Krylov space implemented within the Golub Kahan Bidiagonalization algorithm, here denoted EGKB. RSVD and EGKB are used to efficiently find a low-rank approximation to the SVD of $\mathcal{R}(A)$.

When $A$ is large a limiting factor of the approach to find the KP sum approximation is the need to form the SVD for $\mathcal{R}(A)$ which has the same number of elements as $A$. Even if we decide that we do not need $A$ for subsequent calculations, and immediately overwrite $A$ with $\mathcal{R}(A)$, we are still left with forming an approximate SVD for a large matrix. Thus, for our implementation of the KP sum algorithm using approximate SVDs we also propose replacing  $A$ by its single precision (SP) approximation, hence  reducing the storage by a factor of $2$. More importantly, this allows us to calculate  the approximate SVD for the rearrangement of $A$ in single precision (SP). The SP computations can be performed approximately twice as fast as for double precision objects because SP requires half the storage and uses half the memory compared to the double precision format \cite{higham2022mixed}. Thus, to further speed up the KP sum approximation process, we use SP computations in EGKB and RSVD, denoted EGKB(SP) and RSVD(SP), respectively,  using $32$ bit SP numbers to store the matrix $\mathcal{R}(A)$ and to perform the operations with $\mathcal{R}(A)$.    After obtaining the KP approximation to $A$ with $k$ terms, denoted by $\tilde{A}_k$, with  $k\le R$,  the Kronecker rank $R$ of $A$, the solution of \cref{eq:Tv} is found using the SB algorithm with all uses of $A$ replaced by $\tilde{A}_k$. The conjugate gradient least squares (CGLS) algorithm is applied to solve the required least squares problem within the SB iteration  \cite{paige1982lsqr}. EGKB(SP) and RSVD(SP) are contrasted with the use of the double precision  at all stages of the image reconstruction. Overall, this contrasts with the use of a short KP sum approximation used in \cite{garvey2018singular} to find an SVD approximation to $A$ directly, and the use of structure to find a suitable KP approximation, \cite{kamm2000optimal}, for eventual use in a regularization algorithm. 

Throughout we assume that the problem described by \cref{eq:forward problem} is two dimensional. The vectors $\bfx$ and $\bfb$ represent the two dimensional image and its contaminated version, respectively. To perform matrix vector operations, and hence implement the algorithms efficiently, we will assume when needed the notation that $\bfx = \Vector(\bfX)$ and  $\bfX =\Arrayop(\bfx)$, where $\bfX$ is the image without reshaping to a vector. 
\subsection{Contributions}
We use the EGKB(SP) and RSVD(SP) to efficiently estimate an approximate TSVD for $\mathcal{R}(A)$ rather than using the full SVD of $\mathcal{R}(A)$. In order to show the efficiency of these algorithms, they are contrasted with EGKB and RSVD when applied to the same problem, but with no SP approximations. Empirically, the number of terms $k$ in the Kronecker product sum approximation to $A$ is automatically determined using the decay of the geometric mean of the products of the elements of the bidiagonal matrix in the EGKB(SP). The approximate $\tilde{A}_k$ for $A$  is then used  to solve the total variation regularization problem without needing further use of the matrix $A$. Given the KP approximation, the regularized solution at each outer iteration of the total variation algorithm is computed in double precision using CGLS.  Our numerical results confirm that it is more efficient to use $\tilde{A}_k$  than $A$ in the SB algorithm, and  computational complexity estimates are provided to support this observation. The results confirm the effectiveness of the SP algorithms for obtaining $\tilde{A}_k$, and of the use of the automatic stopping rule for determining $k$ using the EGKB(SP) algorithm.

\subsection{Organization}
This paper is organized as follows. The KP sum approximation to the reordered matrix $\mathcal{R}(A)$ for the non-separable blurring operator $A$ is described in \cref{sec:Approx A}. Finding $\mathcal{R}(A)$ from $A$ is presented in \cref{subsec:reorder A}, and then in \cref{subsec:EGKB,subsec:RVD}, we provide the EGKB and RSVD algorithms that are used to form an SVD approximation for $\mathcal{R}(A)$. The basic idea of the SB iterative algorithm for isotropic and anisotropic total variation regularization  for the solution of \cref{eq:Tv} is reviewed in \cref{subsec:SB,sec:Anisotropic,sec:isotropic}. In \cref{subsec:CGLS}, we describe the CGLS algorithm that is used to solve the outer generalized Tikhonov update  within the SB iteration. 
Comparative computational costs of the EGKB and RSVD algorithms, as well as directly using  matrix $A$  in the CGLS, are provided in \crefrange{subsec:KPapprox}{subsec:CostSB}. In \cref{seubsec:SP}, we discuss using SP for the computations in EGKB and RSVD, and hence for obtaining the SP approximation $\tilde{A}_k$ to $A$.
Numerical results and conclusions are presented in \cref{sec:Numerical examples,sec:Conclusion}, respectively.

\section{The Kronecker Product Approximation}\label{sec:Approx A}
In this section, we discuss using a well-known method to approximate the large matrix $A$ by a sum of Kronecker Products. This approach requires calculating the SVD of the reordering of the entries of $A$ \cite{golub2013matrix}. The resulting algorithm is efficient in terms of storage and time complexity.

\subsection{Rearrangement \texorpdfstring{$\mathcal{R}(A)$ of $A$}{R(A) of A}} \label{subsec:reorder A}
Under the assumption that the operator $A \in \mathbb{R}^{M\times N}$ with $M=m_1 m_2$ and $N=n_1n_2$ is large and not already structured, $A$ can be approximated using a sum of Kronecker products as first proposed in \cite{loan1992approximation}. To find the approximation we must  rearrange the elements of the matrix $A$ into another matrix, namely $\mathcal{R}(A)$. Here we generate $\mathcal{R}(A)$ using the method discussed in \cite{loan1992approximation} and subsequently in \cite[eq.12.3.16]{golub2013matrix}.

Let $A$ form a block partition
\begin{equation}\label{eq:block A for KP}
    A= \begin{bmatrix}
        A_{1,1} & A_{1,2} & \dots & A_{1,n_1}\\
        A_{2,1}& A_{2,2} & \dots & A_{2,n_1}\\
        \vdots&\vdots&\ddots&\vdots\\
        A_{m_1,1}&A_{m_1,2}&\dots&A_{m_1,n_1}
    \end{bmatrix}, \quad A_{i,j} \in \mathbb{R}^{m_2\times n_2}.
\end{equation}
With $\Vector(A)=\begin{bmatrix}
    A^\top_{1,1} & A^\top_{2,1}&\dots&A^\top_{1,2}&\dots&A^\top_{m_1,n_1}
\end{bmatrix}^\top$, one can construct the rearrangement of $A$ as
\begin{equation}\label{eq:nitial rearrange A}
 \mathcal{R}(A)= \begin{bmatrix}
      A_1\\
      A_2\\
      \vdots\\
      A_{n_1}
  \end{bmatrix}, \quad
  A_j= \begin{bmatrix}
      \Vector(A_{1,j})^\top\\
      \Vector(A_{2,j})^\top\\
      \vdots\\
      \Vector(A_{m_1,j})^\top\\
  \end{bmatrix}, \quad j=1:n_1.
\end{equation}
It should be noted that $\mathcal{R}(A)$ has the same number of entries as $A$ but is of size $m_1 n_1 \times m_2 n_2$, which differs from the size of $A$.
\begin{theorem}\cite[Theorem 12.3.1]{golub2013matrix}\label{eq:Kronecker Product SVD}
Let $A \in \mathbb{R}^{M \times N}$, with $M=m_1m_2$ and $N=n_1 n_2$, be blocked as in \cref{eq:block A for KP}, and the SVD of $\mathcal{R}(A)$ be given by 
\begin{equation}\label{eq:SVDofR(A)}
   \mathcal{R}(A)= \widetilde{U}\widetilde{\Sigma}\widetilde{V}^{\top},
\end{equation}  then 
\begin{equation}\label{eq:KPsum}
    A= \sum^R_{i=1} \widetilde{\sigma}_i \  \widetilde{U}_i \otimes \widetilde{V}_i,
\end{equation}
    where $\widetilde{\sigma}_i$ are the singular values of $ \mathcal{R}(A)$, $\widetilde{U}_i = \Arrayop(\widetilde{\bfu}_i, m_1, n_1)$ and $\widetilde{V}_i =\Arrayop(\widetilde{\bfv}_i, m_2, n_2)$\footnote{\text{Notice we specify here that the arrays are reshaped to sizes} $m_1\times n_1$ and $m_2\times n_2$, respectively.} and $R$ is the Kronecker rank of $A$.
\end{theorem}
Let $\tilde{A}_k$ denote the approximation to $A$ using \cref{eq:KPsum} with $k\le R$ terms.  By setting $A_{x_i} = \widetilde{\sigma}_i \  \widetilde{U}_i$ and $A_{y_i} =  \widetilde{V}_i$, we can  write the statement in \Cref{eq:Kronecker Product SVD} as
\begin{equation}\label{eq:KP approx}
     \tilde{A}_R= \sum^R_{i=1} A_{x_i} \otimes A_{y_i},
\end{equation}
where $A_{x_i}\in \mathbb{R}^{m_1 \times n_1}$  and $A_{y_i}\in \mathbb{R}^{m_2 \times n_2} $.

Note that $A$ is large, and accordingly so is $\mathcal{R}(A)$. Thus, finding the full SVD of $\mathcal{R}(A)$ can be computationally unattractive. Instead, we present two efficient known approaches for estimating an approximation of the SVD of any given matrix. 
\subsection{Enlarged-GKB} \label{subsec:EGKB}
The GKB algorithm, also known as the Lanczos bidiagonalization algorithm, was introduced in \cite{golub1965calculating} for reducing any general matrix to bidiagonal form. This method is, historically, one of the most used methods for estimating the SVD of a large matrix. It is well known that the GKB approximation for the SVD yields a spectrum that approximates the total spread of the true singular spectrum, inheriting the condition of the problem and only giving a portion of the dominant spectrum \cite{paige1982lsqr}. Specifically, suppose that $\tilde{\sigma}_i$ is the $i^{\text{th}}$ singular value obtained by the GKB algorithm, for $i=1$ to $k+p=k_p$, then $\tilde{\sigma}_i$ is only a good approximation to the true $\sigma_i$ for some $i=1$ to $k$, where $k \ll k_p$ \cite{paige1982lsqr}. Here $k$ is a target rank and $p$ is a sampling parameter. It was also suggested in \cite[section 3.3.1]{larsen1998lanczos} that using larger $k$ leads to significant improvement in the accuracy of the Ritz values, meaning that larger $k$ yields an accurate approximation to the eigenvalues of the desired matrix. While the approach of  enlarging the projected system was suggested in \cite{larsen1998lanczos}, and  adopted in  \cite{renaut2020fast,renaut2017hybrid} in the context of selecting the regularization parameter and  improving the quality of reconstruction of the geophysics problems, respectively, there does not appear to be explicit wide spread use of the idea to enlarge the Krylov subspace with subsequent truncation in the context of SVD approximation as presented here. We use EGKB to denote the TSVD obtained using the enlarged GKB algorithm.  The steps of the EGKB algorithm are described in \Cref{alg:EGKB}. 
\begin{algorithm}[htb!]
\caption{\textbf{EGKB:  providing a TSVD approximation with at most $k$ terms to a given matrix  $B$ \cite{golub1965calculating}.
}}\label{alg:EGKB}
\begin{algorithmic}[1]
    \Require $B \in \mathbb{R}^{\bar{M}\times \bar{N}}$ with $\bar{M}=m_1 n_1$ and  $\bar{N}= m_2 n_2$; initial guess $\bfy \in \mathbb{R}^{\bar{M}}$; a target rank $k \ll \text{min} (\bar{M},\bar{N})$; a sampling parameter $p$;  oversampled projected problem $k_p=k+p$; $\tau_{EGKB}$
    \State {\text{Initialize:} $t_1=\|\bfy\|_2$, $\bfs_1=\bfy/t_1$,  $\bfq=B^\top \bfs_1$, $\alpha_1=\|\bfq\|_2$, $\bfq_1=\bfq/\alpha_1$, $\nu_2=100$, $\zeta_1=1$, $j=2$ \;}
      \While {not converged ($\nu_j>\tau_{EGKB}$, $\nu_{j+1}<\nu_j$ or $j<k+p+2$) } 
\State {$\bfs_j=B\bfq_{j-1}- \alpha_{j-1} \bfs_{j-1}$ \;}
\State {$t_j=\|\bfs_j \|_2$, \ $\bfs_j=\bfs_j /t_j $ \;}
\State {$\bfq_j=B^\top \bfs_j- t_j \bfq_{j-1}$ \;}
 \State {$\bfq_j = \bfq_j - \sum^{j-1}_{i=1} \langle \bfq_j , \bfq_i\rangle\bfq_i$ \ \ \text{one-sided reorthogonalization} \;}
\State {$\alpha_j=\|\bfq_j \|_2, \ \bfq_j=\bfq_j/\alpha_j$ \;}
\State {$\zeta_j= \alpha_j t_j$,  $\nu_{j+1}=\zeta_j\zeta_{j-1}$}
 \State {$\nu_{j+1} > \nu_j \text{ or }  \nu_j<\tau_{EGKB}$, and $j-1\ne k_p$, set $k=j$.}
\State {$j=j+1$ \;}
\EndWhile
\State{Define $S_{\bar{M}\times {(k_p+1)}}=[\bfs_1,\dots,\bfs_{(k_p+1)}]$ and $ Q_{\bar{N}\times {k_p}}=[\bfq_1,\dots,\bfq_{k_p}]$\;}
\State{Scalars $t_j$ and $\alpha_j$ define a lower
bidiagonal matrix $C_{(k_p+1)\times k_p}$  \;}
\State{Obtain $B = S_{\bar{M} \times (k_p+1)} C_{(k_p+1)\times k_p} Q_{k_p\times \bar{N}}^\top$ \:}
\State {SVD for bidiagonal matrix $C_{(k_p+1) \times k_p}$: \   $U_{(k_p+1)\times k_p} \bar{\Sigma}_{k_p\times k_p} V_{k_p\times k_p} = \texttt{svds}(C_{(k_p+1)\times k_p},k_p)$ \;}
\State {Truncate  $S_{\bar{M} \times (k_p+1)}$, $ U_{(k_p+1)\times k_p}$, $\bar{\Sigma}_{k_p\times k_p}$, $V_{k_p\times k_p}$, and $Q_{\bar{N}\times k_p} $ to the first $k$ terms \;}
\Ensure $\bar{B}_k = ( S_{\bar{M} \times k} U_{k\times k} ) \bar{\Sigma}_{k\times k} (Q_{\bar{N}\times k} V_{k\times k})^\top =\bar{U}_{\bar{M}\times k} \bar{\Sigma}_{k\times k} \bar{V}^\top_{k\times \bar{N}}$
\end{algorithmic}
\end{algorithm}
Note that $p$ is a user-defined parameter, typically taken to be small but increasing with  $k$. To maintain the orthogonality of the Lanczos columns of $Q$, one-sided reorthogonalization is applied at step $6$ as proposed in \cite{simon2000low} using  the modified Gram Schmidt reorthogonalization for the right Lanczos columns defining $Q_{k+1}$. \Cref{alg:EGKB} can easily be modified to also require that step $9$ is only implemented if $j>k_{min}$, where $k_{min}$ is an input parameter that determines a minimum size for $k$ in $\bar{B}_k$.

When we apply $k_p=k+p$ steps of the EGKB to $B=\mathcal{R}(A)$, we find the $k+p$-term  TSVD approximation to $\mathcal{R}(A)$ that is generated for the Krylov subspace of size $k+p$. The  $k$-term approximation, denoted by $\bar{\mathcal{R}}_k(A)$, is obtained by truncation of the $k+p$-term TSVD to $k$-terms. To find a suitable $k_p$ we might calculate the singular values for selected choices of increasing $k_p$, to determine how the singular values are converging. This would be inefficient, since it would require multiple SVD calculations for selected choices of $j$. Instead we introduce a new stopping criteria 
using the decay of the products 
$\zeta_j= \alpha_j t_j$. These decay to zero, but not monotonically,  with $j$ \cite[Corollary 3.1]{gazzola2016lanczos}. Therefore, we use the decay of the geometric mean of $\zeta_j$ calculated only over two terms, using  
$\nu_j=\sqrt{\zeta_j \zeta_{j-1}}$, defined for $j>2$. We estimate $k$ in terms of the minimum $j$ that satisfies 
\begin{equation}\label{eq:stopnu}
k=\min_j \{\nu_{j+1} > \nu_j \text{ or }  \nu_j<\tau_{EGKB}\}.
\end{equation}
This determines the number of terms $k$ to use for $\bar{\mathcal{R}}_k(A)$. After $k$ is chosen in step $9$ the EGKB algorithm is not terminated, but it continues to run up to $k_p=k+p$ terms, for the adjusted $k$. At step $16$  the obtained SVD is truncated to $k$ terms, yielding the desired $\bar{\mathcal{R}}_k(A)$.

\subsection{RSVD} \label{subsec:RVD}
The RSVD algorithm \cite{halko2011finding} uses random sampling to build a low dimensional subspace to capture the majority of the action of the matrix $B=\mathcal{R}(A)$. Then $B$ is  restricted to approximation in a smaller subspace that is more computationally amenable to decomposition using the SVD. The main steps of the RSVD algorithm for a general matrix $B$ are shown in \Cref{alg:RSVD}.
\begin{algorithm}[htb!]
\caption{\textbf{RSVD: providing a  TSVD approximation with $k$ terms to a given matrix  $B$ of size $\bar{M}$ by  $\bar{N}$, where $\bar{M}\ge \bar{N}$} \cite{halko2011finding}.\label{alg:RSVD}}
\begin{algorithmic}[1]
\Require $B\in \mathbb{R}^{\bar{M}\times \bar{N}}$ with $\bar{M}=m_1 n_1$ and  $\bar{N}= m_2 n_2$; a sampling parameter $p$; a target rank $k \ll \text{min} (\bar{M},\bar{N})$, oversampled projected problem $k_p=k+p$  
\State\text{Draw a random Gaussian matrix $F \in \mathbb{R}^{\bar{N}\times k_p}$\;}
\State\label{step2:RSVD}\text{Sample matrix $Y=(B B ^\top)^qB F \in \mathbb{R}^{\bar{M}\times k_p}$ (implemented  with orthonormalization)\;}
\State\label{step3:RSVD}\text{Construct $Q \in \mathbb{R}^{\bar{M}\times k_p}$ whose columns form an orthonormal basis for the range of $Y$ \;}
\State\label{step4:RSVD}\text{Form the smaller matrix $H = Q^\top B \in \mathbb{R}^{k_p \times \bar{N}}$ \;}
\State\label{step5:RSVD}\text{Compute the SVD of $H$, $U \in \mathbb{R}^{k_p \times k_p}$, $\Sigma \in \mathbb{R}^{k_p \times {k_p}}$, and $V \in \mathbb{R}^{\bar{N} \times {k_p}}$ \;}
\State\label{step6:RSVD}\text{Form the matrix $ \bar{U} = Q  U \in \mathbb{R}^{\bar{M} \times k_p}$ \;}
\State\label{step8:RSVD}{Truncate  $\bar{U}_{\bar{M} \times k_p}$, $\Sigma_{k_p\times k_p}$, and $V_{ \bar{N}\times k_p}$ to the first $k$ terms   \;}
\Ensure $\bar{B}_k = \bar{U}_{\bar{M} \times k} \bar{\Sigma}_{k\times k} \bar{V}_{k\times \bar{N}}^\top$
\end{algorithmic}
\end{algorithm}

The matrix $F$ reduces the size of $B$ in step $1$ of \Cref{alg:RSVD} and $Y$ in step $2$ contains as much information of $B$ as possible. One of the features of this algorithm is that the target rank $k$ and sampling parameter $p$ can be specified in advance. \Cref{alg:RSVD} expands the size of the projected problem by using the larger quantity $k_p=k+p$ rather than just $k$. Thus, the quality of $\bar{\mathcal{R}}_k(A)$ as an approximation to $B=\mathcal{R}(A)$ depends critically on the choices of $k$ and $p$. The RSVD in this work uses the estimate $k$ obtained from the EGKB with the new stopping criteria \cref{eq:stopnu}. The extra $p$ terms are removed from the SVD approximation, as given in step~$7$ of \Cref{alg:RSVD}, yielding the $k$-term RSVD approximation,  $\bar{\mathcal{R}}_k(A)$.

Practically, to reduce computational error at step $2$ of \Cref{alg:RSVD}, potentially leading to contamination of the basis $Q$ \cite{halko2011finding}, step $2$ is implemented  with orthonormalization of the columns of $Y$ at every single application of $B$ and $B^\top$. 

\section{Isotropic and Anisotropic Regularization}\label{sec:TVReg}
We briefly present the  isotropic and anisotropic total variation regularization algorithms  based on the Split Bregman (SB) algorithm of Goldstein and Osher \cite{goldstein2009split}. These algorithms require an outer iteration for each step of the SB algorithm and an inner iteration to solve the associated least squares system of equations that arises each outer iteration. We use superscripts $\ell$ to denote the outer iterations and $\iota$ to denote the inner iterations, as needed. Furthermore, in each case we need to impose a maximum number of outer iterations and a tolerance for convergence of the outer iteration. We use limiting number of iterations $L_{max}$ and $I_{max}$, and convergence tolerances   $\tau_{SB}$ and $\tau_{CGLS}$, for the SB algorithm and CGLS solve, respectively. 
\subsection{The Split Bregman Algorithm}\label{subsec:SB}
The SB approach for finding the solution of \cref{eq:Tv} in one dimension involves introducing $\bfd=L \bfx$, leading to the constrained problem given by 
\begin{equation}\label{eq:constSB}
     \min_{\bfx,\bfd}\{\|A \bfx - \bfb \|^2_2 + \beta \|\bfd\|_1\} \quad \text{such that} \quad \bfd=L \bfx.
\end{equation}
Using the augmented Lagrangian approach to solve \cref{eq:constSB} we  arrive at the  unconstrained optimization problem
\begin{equation}\label{eq:split_Bregman}
   (\bfx,\bfd)= \arg  \min_{\bfx,\bfd}\{\|A \bfx - \bfb \|^2_2 + \lambda^2 \| L \bfx - \bfd \|^2_2 + \beta \|\bfd\|_1 \},
\end{equation}
where $\lambda>0 $ is a regularization parameter. 
Note that due to the splitting of the $\ell_1$ and $\ell_2$ components of this functional, an efficient minimization can be performed iteratively for $\bfx$ and $\bfd$ separately.  As proposed in \cite{goldstein2009split}, the SB iterative algorithm for solving \cref{eq:split_Bregman} with respect to $\bfx$ and $\bfd$ is given by
\begin{align}\label{eq:update x}
      \bfx^{(\ell+1)}&= \arg  \min_{\bfx}\{\|A \bfx - \bfb \|^2_2 + \lambda^2 \|L \bfx-(\bfd^{(\ell)}-\bfg^{(\ell)})\|^2_2 \},\\\label{eq:update d}
      \bfd^{(\ell+1)}&= \arg  \min_{\bfd}\{ \lambda^2 \|\bfd -(L \bfx^{(\ell+1)}+\bfg^{(\ell)})\|^2_2 + \beta \|\bfd\|_1 \},
\end{align}
with the update formula 
\begin{gather}\label{eq:update g}
    \bfg^{(\ell+1)}=\bfg^{(\ell)}+L \bfx^{(\ell+1)}-\bfd^{(\ell+1)},
\end{gather}
 where  $\bfg^{(\ell+1)}$ is generated by the Bregman distance. This defines the solution for $\bfx^{(\ell+1)}$ as an outer iterative solution for iteration $\ell+1$, accompanied by an inner iterative scheme to find $\bfx^{(\ell+1)}$ as the solution of \cref{eq:update x}. Notice that the matrices in the SB iteration as given are not iteration dependent. Therefore, in some cases it is possible to use a joint decomposition of $A$ and $L$ to solve this problem directly at each iteration, \cite{sweeney2024parameterselectiongcvchi2}. A joint decomposition is not the focus of our work, rather for the two dimensional anisotropic and isotropic regularization algorithms presented in \cref{sec:Anisotropic,sec:isotropic}, respectively, in which we replace $A$ by a $\tilde{A}_k$ with $k$ terms, we  use an iterative solver. 

\subsection{Anisotropic TV}\label{sec:Anisotropic}
For the two dimensional image $\bfX$, the regularization can be applied by considering each direction of the image separately. Suppose that $I$ is an appropriately sized identity matrix and that $\otimes$ is used to denote the Kronecker product. Then the \textbf{anisotropic}  TV regularization for $\bfx=\Vector({\bfX})$ is given by 
\begin{equation*}
      \min_\bfx\{\|A \bfx - \bfb \|^2_2 + \beta_x \|( L_x \otimes I) \bfx\|_1 + \beta_y \|(I \otimes L_y) \bfx\|_1\},
\end{equation*}
where $\beta_x$ and $\beta_y$ denote the regularization parameters associated with the derivatives for $\bfx$ in horizontal and vertical directions, $L_x$ and $L_y$, respectively. Replacing $(L_x \otimes I) \bfx$ and $(I \otimes L_y) \bfx$ by $\bfd_x$ and $\bfd_y$ respectively, leads to the constrained representation 
\begin{equation*}
      \min_\bfx\{\|A \bfx - \bfb \|^2_2 + \beta_x \|\bfd_x\|_1 + \beta_y \|\bfd_y \|_1\}  \quad \text{such that} \quad \bfd_x = (L_x \otimes I) \bfx \quad \text{and} \quad \bfd_y =(I\otimes L_y) \bfx.
\end{equation*}
Using the regularization parameters $\lambda_x$ and $\lambda_y$ for each constraint, respectively, we can write the complete minimization problem as 
\begin{equation}\label{eq:anispt_with SB iter agu}
     \min_{\bfx,\bfd_x,\bfd_y}\{\|A \bfx - \bfb \|^2_2 +  \lambda^2_x\| (L_x\otimes I) \bfx -  \bfd_x \|^2_2 + \lambda^2_y\| (I \otimes L_y) \bfx - \bfd_y  \|^2_2 +  \beta_x \| \bfd_x \|_1 +\beta_y \| \bfd_y \|_1 \}.
\end{equation}

The SB iteration as given in \crefrange{eq:update x}{eq:update g} was introduced for the solution of minimization problem \cref{eq:anispt_with SB iter agu}  in \cite{goldstein2009split}. Specifically, the solution of \cref{eq:anispt_with SB iter agu} is found by iterating over the following subproblems 
\begin{align}
 \begin{split}
     \bfx^{(\ell+1)}=& \argmin_{\bfx}  \{\|A \bfx - \bfb \|^2_2 + \lambda^2_x \|  (L_x\otimes I) \bfx - (\bfd_x^{(\ell)} - \bfg_x^{(\ell)})  \|^2_2  \\&\qquad \qquad \qquad \qquad \  + \lambda^2_y \|  (I\otimes L_y) \bfx - (\bfd_y^{(\ell)} - \bfg_y^{(\ell)})  \|^2_2\},\label{eq:x in 2d}
     \end{split}
     \\ 
     \bfd_x^{(\ell+1)}=&   \argmin_{\bfd_x} \{ \lambda^2_x \| \bfd_x - \bfc_x^{(\ell)} \|^2_2 + \beta_x \| \bfd_x\|_1 \}\nonumber,
          \\
     \bfd_y^{(\ell+1)}=&   \argmin_{\bfd_y} \{ \lambda^2_y \| \bfd_y - \bfc_y^{(\ell)} \|^2_2 + \beta_y \| \bfd_y\|_1 \}\nonumber, 
\end{align}
where 
    $\bfc^{(\ell)}_x=(L_x\otimes I) \bfx^{(\ell+1)}+\bfg_x^{(\ell)}$ 
and $    \bfc^{(\ell)}_y=(I\otimes L_y) \bfx^{(\ell+1)}+\bfg_y^{(\ell)}$.  
Each iteration is completed with the update for $\bfg$ given by
\begin{align}\label{eq:update g_x ani}
    \bfg^{(\ell+1)}_x&=  \bfg^{(\ell)}_x + (L_x\otimes I) \bfx^{(\ell+1)} -\bfd^{(\ell+1)}_x\\\label{eq:update g_y ani}
    \bfg^{(\ell+1)}_y&=  \bfg^{(\ell)}_y + (I\otimes L_y) \bfx^{(\ell+1)} -\bfd^{(\ell+1)}_y, 
\end{align}  
which is equivalent to \cref{eq:update g} but explicitly framed for the two dimensional formulation.

Now, by defining 
$
    \gamma_x=\lambda^2_x/\beta_x$,  and   $\gamma_y=\lambda^2_y/\beta_y$, 
we can write the update for $\bfd_x$ and $\bfd_y$ as 
$     \bfd^{(\ell+1)}_x= \argmin_{\bfd_x} \{  \gamma_x\|\bfd_x - \bfc^{(\ell)}_x) \|^2_2 +  \|\bfd_x\|_1 \}$ and  
          $\bfd^{(\ell+1)}_y=   \argmin_{\bfd_y} \{  \gamma_y\|\bfd_y - \bfc^{(\ell)}_y) \|^2_2 +  \|\bfd_y\|_1 \}$. 
To find each element in $\bfd_x$ and $\bfd_y$ we use the general shrinkage formula
\begin{equation}\label{eq:d thresholding ani}
     (\bfd^{(\ell+1)})_i= \text{shrink}\left((\bfc^{(\ell)})_i,\frac{1}{\gamma} \right),
\end{equation}
where $\text{shrink}(w,v) =\text{sign}(w) \, \text{max} (|w|-v,0)$, and \cref{eq:d thresholding ani} solves the general problem 
\begin{align*}
     \bfd^{(\ell+1)}=&  \argmin_{\bfd} \{  \gamma\|\bfd - \bfc^{(\ell)}) \|^2_2 +  \|\bfd\|_1 \}.
\end{align*}

In summary, the SB algorithm for  anisotropic TV regularization is given in  \Cref{alg:Anisotropic}. 

\begin{algorithm}[htb!]
\caption{\cite{goldstein2009split} \textbf{Split Bregman Anisotropic Total Variation}}\label{alg:Anisotropic}
\begin{algorithmic}[1]
\Require $A$; $\bfb$; $L_x$; $L_y$ $ \lambda_x$; $ \lambda_y$; $\beta_x$; $\beta_y$; $L_{max}$; $\tau_{SB}$
    \State {\text{Initialize:} $\ell=0, \ \bfx^{(\ell)} = \bfb , \ \bfd_x^{(\ell)}= \bfd_y^{(\ell)} = \bfg_x^{(\ell)} = \bfg_y^{(\ell)} = 0 $. \;}
\While {not converged and $\ell<L_{max}$}   
\State\label{alg:WhileAni1} \text{Update $\bfx$: solve   \cref{eq:x in 2d} \;}
\State\label{alg:WhileAni2} \text{Update $\bfd_x$ and $\bfd_y$: using    \cref{eq:d thresholding ani} \;}
\State\label{alg:WhilAni4} \text{Update $\bfg_x$ and $\bfg_y$:   using  \cref{eq:update g_x ani} and \cref{eq:update g_y ani}, resp.\;}
\State \text{Increment $\ell=\ell+1$\;}
\EndWhile
\Ensure {$\bfx$ \;} 
\end{algorithmic}
\end{algorithm}

\subsection{Isotropic TV}\label{sec:isotropic}
The isotropic TV formulation for the regularized solution of \cref{eq:forward problem} is given by
\begin{equation*}
      \min_\bfx\{ \|A \bfx - \bfb \|^2_2 +\beta  \|\sqrt{|(L_x \otimes I)\bfx|^2 + |(I \otimes L_y)\bfx|^2} \|_1\}.
\end{equation*}
As for the anisotropic case, we introduce $\bfd_x$ and $\bfd_y$ to satisfy the constraints   $
 \bfd_x = (L_x \otimes I) \bfx$  and  $\bfd_y =(I\otimes L_y) \bfx$, which then leads to the unconstrained minimization problem
\begin{equation}\label{eq:ini_isto}
     \min_{\bfx,\bfd_x,\bfd_y}\{\|A \bfx - \bfb \|^2_2 + \lambda^2_x \| (L_x \otimes I)\bfx - \bfd_x   \|^2_2 + \lambda^2_y \| (I \otimes L_y) \bfx-\bfd_y  \|^2_2 +\beta \|\sqrt{|\bfd_x|^2 + |\bfd_y|^2 }\|_1 \}.
\end{equation}
Now the SB formulation is modified as compared to the anisotropic case through the updates that are required for $\bfd_x$ and $\bfd_y$.  Here, $\bfd_x$ and $\bfd_y$  do not decouple and there is a single parameter $\beta$, while the updates for $\bfx$ and the Lagrange multipliers are unchanged from the anisotropic formulation. Thus, we consider the problem \cref{eq:x in 2d} with the following subproblem
\begin{equation}\label{eq:dx dy istrophic}
\begin{split}
 (\bfd^{(\ell+1)}_x,\bfd^{(\ell+1)}_y) =&   \argmin_{\bfd_x,\bfd_y} \{ \lambda^2_x \| \bfd_x - \bfc^{(\ell)}_x \|^2_2 + \lambda^2_y \| \bfd_y - \bfc^{(\ell)}_y \|^2_2 
 +\beta \|\sqrt{|\bfd_x|^2 + |\bfd_y|^2 }\|_1 \}\\
=& \argmin_{\bfd_x,\bfd_y} \{ \tilde{\gamma}_x \| \bfd_x - \bfc^{(\ell)}_x \|^2_2 + \tilde{\gamma}_y \| \bfd_y - \bfc^{(\ell)}_y \|^2_2 
 + \|\sqrt{|\bfd_x|^2 + |\bfd_y|^2 }\|_1 \}, 
\end{split}
\end{equation}
where $\tilde{\gamma}_x = \lambda^2_x/\beta$ and $\tilde{\gamma}_y = \lambda^2_y/\beta$. 
The solutions of the minimization problem \cref{eq:dx dy istrophic} can be expressed using the shrinkage expressions \cite[eqs 4.4-4.6]{goldstein2009split} 
\begin{equation}\label{eq:d thresholding iso}
     (\bfd^{(\ell+1)})_i= \frac{(\bfc^{(\ell)})_i}{(\bfs^{(\ell)})_i} \, \text{max}\left((\bfs^{(\ell)})_i-\frac{1}{\tilde{\gamma}},0\right) ,
     \end{equation}
where $(\bfs^{(\ell)})_i = \sqrt{|(\bfc^{(\ell)}_x)_i|^2 + | (\bfc^{(\ell)}_y)_i|^2}$.

In summary, the SB algorithm for the isotropic TV regularization is given in  \Cref{alg:Isotropic}.

\begin{algorithm}[htb!]
\caption{\cite{goldstein2009split} \textbf{Split Bregman Isotropic Total Variation}}\label{alg:Isotropic}
\begin{algorithmic}[1]
\Require $A$; $\bfb$; $L_x$; $L_y$; $ \lambda_x$; $ \lambda_y$; $\beta_x$; $\beta_y$; $L_{max}$;  $\tau_{SB}$
\State {\text{Initialize:} $\ell=0, \ \bfx^{(\ell)} = \bfb , \ \bfd_x^{(\ell)}= \bfd_y^{(\ell)} = \bfg_x^{(\ell)} = \bfg_y^{(\ell)} = 0 $. \;}
\While {not converged and $\ell<L_{max}$}
\State\label{alg:Whileis1} \text{Update $\bfx$: solve    \cref{eq:x in 2d}   \;}
\State\label{alg:Whileis2} \text{Update $\bfd_x$ and $\bfd_y$: using    \cref{eq:d thresholding iso}\;}  
\State\label{alg:Whileis4} \text{Update $\bfg_x$ and $\bfg_y$:   using  \cref{eq:update g_x ani} and \cref{eq:update g_y ani}, resp.\;}
\State \text{Increment $\ell=\ell+1$\;}
\EndWhile
\Ensure {$\bfx$ \;} 
\end{algorithmic}
\end{algorithm}

\subsection{CGLS Algorithm}\label{subsec:CGLS}
The Conjugate Gradients Least Square (CGLS) algorithm, presented in \cite[section 7.1]{paige1982lsqr}, is considered for the update \cref{eq:x in 2d}. This algorithm is derived from the well-known Conjugate Gradients algorithm introduced in \cite{hestenes1952methods} and  is widely used for ill-posed problems due to its numerical stability \cite[Section 8.3]{saad2003iterative}. In practice, CGLS is used to solve the  augmented system  derived from  \cref{eq:x in 2d} given by 
\begin{equation}\label{eq:agum x}
    \bfx^{(\ell+1)}=  \argmin_{\bfx}  \displaystyle\{\| 
    \hat{A} \bfx -  
    \hat{\bfb} \|^2_2\},
\end{equation} 
where 
\begin{align}\label{eq:agu A}
   \hat{A}= & \begin{bmatrix}
    A\\  \lambda_x  (L_x\otimes I)\\ \lambda_y  (I\otimes L_y)
    \end{bmatrix},   \ \ \ 
    \hat{\bfb}= \begin{bmatrix}
        \bfb \\ \lambda_x\bfh^{(\ell)}_x \\ \lambda_y\bfh^{(\ell)}_y
        \end{bmatrix}, 
\end{align}
 with $\bfh^{(\ell)}_x = \bfd_x^{(\ell)} - \bfg_x^{(\ell)}$,  
and $\bfh^{(\ell)}_y =  \bfd_y^{(\ell)} - \bfg_y^{(\ell)}$. 
Here we follow the presentation given in \cite[section 7.1]{paige1982lsqr}, leading to
\Cref{alg:CGLS}.
 \begin{algorithm}[htb!]
 \caption{\textbf{The CGLS algorithm to obtain the solution $\bfx$ of \cref{eq:x in 2d} using \cref{eq:agum x}}\label{alg:CGLS}}
 \begin{algorithmic}[1]
 \Require $\hat{A} \in \mathbb{R}^{T \times N}$, with $T=M+2P$; $\hat{\bfb} \in \mathbb{R}^T$;  $I_{max}$ and $\tau_{CGLS}$
 \State {\text{Initialize:} $\iota=0, \ \bfx^{(\iota)} = 0 , \ \bfr^{(\iota)}= \hat{\bfb}, \ \bfw^{(\iota)} =\bff^{(\iota)} =\hat{A}^\top \bfr^{(\iota)}, \ \tau_\iota=\|\bff^{(\iota)}\|^2_2$ \;}
\While {not converged and $\iota<I_{max}$} 
  \State{$\bfz^{(\iota)} = \hat{A} \bfw^{(\iota)}$}
  \State {$\mu_{\iota}  =\tau_\iota / \| \bfz^{(\iota)} \|^2_2 $\;}
 \State{$\bfx^{(\iota+1)} = \bfx^{(\iota)} + \mu_{\iota} \bfw^{(\iota)}$ \;}
 \State{$\bfr^{(\iota+1)} = \bfr^{(\iota)} -\mu_{\iota} \bfz^{(\iota)}$ \;}
\State{$\bff^{(\iota+1)}=\hat{A}^\top \bfr^{(\iota+1)} $ \;}
\State{$\tau_{\iota+1}=\|\bff^{(\iota+1)}\|^2_2 $ \;}
  \State{$\delta_\iota = \tau_{\iota+1}/ \tau_\iota$ \;}
  \State{$\bfw^{(\iota+1)} = \bff^{(\iota+1)} + \delta_\iota \bfw^{(\iota)} $ \;}
 \State{$\iota=\iota+1$ \;}
 \EndWhile
 \Ensure {$\bfx$ \;} 
 \end{algorithmic}
 \end{algorithm}

As seen in \Cref{alg:CGLS}, an advantage of using CGLS is that no reorthogonlization or explicit storage of vectors is needed. The dominant cost is the matrix-vector operations with $\hat{A}$ or $\hat{A}^\top$ required at steps $3$ and $7$, respectively. Thus, at iteration $\ell$ of the outer iteration the computational cost is dominated by the two matrix vector operations with $\hat{A}$ and $\hat{A}^\top$, each dominated by  $\mathcal{O}(2TN)$ flops, yielding that CGLS requires approximately 
\begin{equation}
    \mathcal{O}(4TN\iota)
\end{equation}
flops, when there are $\iota$ iterations to convergence. 

It is pertinent to observe that the LSQR algorithm, based on the EGKB as given in \Cref{alg:EGKB}, can also be used to solve a system such as \cref{eq:agu A}. In particular, for the LSQR we would potentially apply the EGKB directly to the matrix $B=A$, or its KP approximation with $\tilde{A}_k=B$, to obtain a low rank approximation to the system matrix $A$ in  the fit to data term, as in, for example,  \cite{gazzola2016lanczos,paige1982lsqr,renaut2017hybrid}. This works well for standard Tikhonov regularization for which the norm for the solution on the projected space conserves the norm on the full space. For the generalized Tikhonov problem, as given in \cref{eq:agu A}, various approaches are available to correctly solve with the term $\|L \bfx\|_2^2$, including techniques described carefully in the review paper \cite{ChungGazzolareview}. For the purposes of flexibility with the algorithm in relation to adding the anisotropic and isotropic regularizers, we use CGLS, noting that due to the change in the right hand side in \cref{eq:agu A} at each iteration, we would also need to build the Krylov space for the solution for each $\ell$, unless recycling approaches are applied \cite{baglamareichel,ChanWan,soodhalter2020survey}.

\section{Computational costs} \label{sec:comp costs}
The computational cost of the SB algorithm is dominated by the outer generalized Tikhonov update defined by \cref{eq:agum x} for the augmented matrix $\hat{A}$  given in \cref{eq:agu A} at each step. In our experiments, \cref{eq:agum x} is solved using the CGLS algorithm in which all the operations with the system matrix $A$, and its transpose, are implemented using $\tilde{A}_k$.  Thus, we first 
discuss the computational cost of finding  a truncated SVD approximation  $\bar{B}_k$ with $k$ terms to a general matrix $B$ of size $\bar{M}$ by $\bar{N}$ using  the EGKB and RSVD algorithms. Then, we consider the cost of the CGLS algorithm for the solution of \cref{eq:agum x} using the $\tilde{A}_k$ as compared to the direct use of $A$ and its transpose.

\subsection{Computational Cost of the KP sum approximation}\label{subsec:KPapprox}
 For determining the costs of forming an approximate truncated SVD to a general matrix $B =\mathcal{R}(A)$ of size $\bar{M}$ by $\bar{N}$, we suppose that $B$ has no special structure, and operations with $B$, such as matrix-matrix and matrix-vector multiplications, are implemented in the standard approach.  The dominant costs of EGKB and RSVD are compared for the same choices of $k$ and $p$. In practice,  the EGKB may terminate for a smaller $k$ if the stopping test on $\nu_j$ is satisfied for $j<k$.  

\subsubsection{The EGKB Algorithm}\label{subsubsubsec:EGKB Cost} In \Cref{alg:EGKB} steps  $3$ and $5$  use matrix-vector multiplications requiring  $\mathcal{O}(4\bar{M}\bar{N}k_p)$ flops. At the $j^{\text{th}}$ iteration in step $6$, the column $\bfq_j$ is reorthogonalized against all previous columns, requiring  approximately 
$\mathcal{O} (\sum^{k_p}_{j=1} 4\bar{N} j)  \approx \mathcal{O}(2\bar{N}k_p^2)$ flops.
The costs of the other steps within the while loop are of lower order and are thus negligible. At step $15$  the SVD of the lower bidiagonal matrix $C_{k_p} \in \mathbb{R}^{{(k_p+1)}\times k_p}$ is obtained, requiring $\mathcal{O}(4(k_p+1)k_p^{2} + 8k_p^3) \approx \mathcal{O}( 12k_p^3)$ flops, as given in  \cite[Table 5.5.1]{golub2013matrix}. 

If the matrices of left and right singular vectors are required explicitly, as they are to give $\tilde{A}_k$ the evaluation of the left singular vectors $(\bar{U}_{\bar{M}\times k})_i$ for $i=1, \dots, k$ requires $\mathcal{O}(2\bar{M}k^2)$ flops. Similarly, the right singular matrix  $(\bar{V}_{\bar{N} \times k})_i$  takes approximately $\mathcal{O}(2\bar{N}k^2)$ flops.  
Combining all steps together, the dominant complexity of determining the approximate TSVD with $k$ terms using the EGKB algorithm is given by 
\begin{equation}\label{eq:EGKB cost1}
    \text{TSVD of $B$ with \Cref{alg:EGKB}} = 4\bar{M}\bar{N}k_p + 2\bar{M}k^2 + 2\bar{N}(k^2 + k_p^2) + 12k_p^3 = \mathcal{O}(4\bar{M}\bar{N}k_p). 
\end{equation}

\subsubsection{The RSVD Algorithm}\label{subsubsubsec:RSVD Cost} For \Cref{alg:RSVD} we consider the computational cost assuming just one  power iteration ($q=1$). Step $3$, which is written as $(BB^\top) B F$, is implemented efficiently with reorthogonalization at each intermediary  step for a matrix $Y$, as recommended in \cite{halko2011finding}.  Two multiplications with $B$ and one with $B^\top$ cost about $\mathcal{O}(6\bar{M}\bar{N}k_p)$, while  quoting from \cite[Table 5.5.1]{golub2013matrix}, orthogonalizing the columns of the resulting three matrices costs approximately  $\mathcal{O}(9\bar{M}k_p^2-3k_p^3)$ flops. Forming $H=Q^\top B$ at step $4$ adds another  $\mathcal{O}(2\bar{M}\bar{N}k_p)$ flops, and determining 
the economic SVD at step $5$ is an additional $\mathcal{O}(\bar{N}{k^2_p})$ flops. Since we seek a truncated left singular matrix $\bar{U}_{\bar{M} \times k}$  in step $7$ of \Cref{alg:RSVD}, we only need the first $k$ columns of $\bar{U}$ in step $6$. Practically, when we do not need to calculate $\bar{U}$,  we can store both $Q$ and the first $k$ columns of $U$  in order to implement forward and transpose products  with  $\bar{U}$. If we need the matrix directly, an additional $\mathcal{O}(2\bar{M}k_pk)$ flops are required. 

Combining all steps together, the dominant complexity of determining the approximate TSVD with $k$ terms using the RSVD algorithm is given by  
\begin{equation}\label{eq:RSVD cost}
    \text{TSVD of $B$ with \Cref{alg:RSVD}} = 8\bar{M}\bar{N}k_p + 2\bar{M}k_pk + \bar{N} k_p^2 -3k_p^2 + 9\bar{M} k_p^2= \mathcal{O}(8\bar{M}\bar{N}k_p).  
\end{equation}

\subsubsection{The Cost of forming the terms in KP}\label{subsubsubsec:terms in KP Cost}
Assume the truncated SVD to $\mathcal{R}(A)$ is available, then we need to look at the computational cost for forming  $\tilde{A}_k$, as given in \cref{eq:KPAx}. To form the $i^{\text{th}}$ term in the KP $A_{x_i}$, we need to calculate the product $\bar{\sigma_i} \bar{U}_i\approx \widetilde{\sigma}_i \widetilde{U}_i$. Because $\bar{U}_i$ is size  $m_1 \times n_1$, this costs $\mathcal{O}(m_1 n_1) = \mathcal{O}(\bar{M})$ flops. Note that the term $A_{y_i}$ is only a reshaping of $\bar{V}_i$ and hence does not require a significant computational effort. Since there are $k$ terms in $\tilde{A}_k$,  the total cost is $\mathcal{O}(\bar{M} k)$ flops.

It is evident from the discussion in \crefrange{subsubsubsec:EGKB Cost}{subsubsubsec:terms in KP Cost} that the cost of finding the TSVD approximation to $\mathcal{R}(A)$ dominates the cost of forming $\tilde{A}_k$.

\subsection{Computational Cost of the CGLS algorithm}\label{subsec:CGLScost}
The results of the EGKB and RSVD are used to produce the alternative KP summations, $\tilde{A}_k$,   to $A$ with $k$ terms instead of the full SVD  which yields a summation with a full Kronecker rank ($k=R$). Thus, computational savings are also achieved using these algorithms to implement the matrix-vector multiplications required in the CGLS algorithm for the inner iteration. To investigate the costs of the operations with a $\tilde{A}_k$, 
we consider a dense matrix $A \in \mathbb{R}^{M \times N}$,  
and an array $\bfX$  of size $n_2\times n_1$, such that $\Arrayop(\bfx)=\bfX$.  Then 
  \begin{equation}\label{eq:KPAx}
     \Arrayop( A\bfx)  = \sum^R_{i=1} A_{y_i}\bfX A^\top_{x_i} \in \mathbb{R}^{m_2 \times m_1}.
       \end{equation} 
       
The operations \cref{eq:KPAx} require $\mathcal{O}\left((2 m_2n_1n_2+2 m_1m_2n_1)R\right) = \mathcal{O}\left(2( m_2 N+ Mn_1)R\right) $ flops. Moreover, the operations involving $A^\top$ have an equivalent computational cost. Using the EGKB, or RSVD for $B=\mathcal{R}(A)$ reduces this cost to about $\mathcal{O}\left(2( m_2 N+ Mn_1)k\right) $ flops, where we assume $k\le R$. Observing that the regularization matrices are mostly sparse, we ignore their contributions to the forward or adjoint system matrix operations that are required in \Cref{alg:CGLS}.  
Equivalently, we assume the dominant costs for operations with $\hat{A}$ and $\hat{A}^\top$ as defined in \cref{eq:agu A} are due to the operations with $A$ and $A^\top$, or respectively with $\tilde{A}_k$ and $\tilde{A}_k^\top$.   

\subsection{Total Cost of the SB algorithm}\label{subsec:CostSB}

We suppose that  the SB algorithm requires $\ell_{\text{end}}$ steps in total, and that the total number of inner iterations  of the CGLS  algorithm is given  by $\iota_{total}=\sum_{\ell=1}^{\ell_{\text{end}}} \iota_\ell$, where $\iota_\ell$ is the number of iterations taken by CGLS at SB iteration $\ell$. Then the total cost  to convergence depends on the costs of doing matrix operations with the augmented matrix $\hat{A}$, or its approximation in the case of the RSVD and EGKB algorithms, where we also assume that costs of the update steps \cref{eq:d thresholding ani,eq:d thresholding iso} for the anisotropic and isotropic cases, respectively, are negligible.

We present the results for three cases, (i) solving the system directly without forming $\mathcal{R}(A)$, 
(ii) finding the KP approximation with $k$ terms by using the EGKB algorithm with extension of the space to size $k_p=k+p$ to find the approximate SVD, and (iii)  finding the KP approximation with $k$ terms by using the RSVD algorithm to approximate the SVD, using an oversampling parameter $p$ and one power iteration. For \Cref{Tab:costs}, we assume 
$\bar{M}\approx \bar{N}$, and that $p< k $, such that $k\approx k_p\ll \bar{N}$. Therefore, the dominant cost to complete the SB iterations in all cases is the sum of the three columns in \Cref{Tab:costs}.

\begin{table}[ht!]
\footnotesize
\caption{\label{Tab:costs} The costs of reconstruction obtained with $\ell_{\text{end}}$ iterations of the SB algorithm when using (i) $A$ directly, 
(ii) when forming the rank $k$ approximation of $\mathcal{R}(A)$ obtained using the EGKB approximation with size $k+p$ and (iii)  when forming the rank $k$ approximation of $\mathcal{R}(A)$ obtained using the RSVD approximation obtained with oversampling parameter $p$. Recall that the system in \cref{eq:agum x} is of size $T\times N$, where $T$ accounts for the augmented matrix $\hat{A}$ with $T=M +2P$ rows, for the regularization matrix of size $P \times N$.}
\begin{center}
\begin{tabular}[t]{c c c c c}
\toprule
{Method}& &TSVD of $B=\mathcal{R}(A)$  &Form $\tilde{A}_k$&Total CGLS Iterations  \\ \midrule 
{(i) Direct}& & & &$4 (TN)\iota_{total}$ \\ 
 {(ii)  EGKB rank $k$}& &  
 $4\bar{N}^2k_p$ &$\mathcal{O}(\bar{M} k)$ &$4(( m_2 N+ Mn_1)k)\iota_{total}$  \\
(iii) RSVD rank $k$& &  
$8\bar{N}^2k_p$&$\mathcal{O}(\bar{M} k)$& $4(( m_2 N+ Mn_1)k)\iota_{total}$  \\ \bottomrule 
 \end{tabular}
\end{center}
\end{table}

As shown in \Cref{Tab:costs}, the EGKB and RSVD are efficient for the SVD of $B$ when $k$ is small relatively to $\bar{N}$, and  the cost of the RSVD is approximately double that of EGKB for estimating the rank $k$ TSVD of $\mathcal{R}(A)$, see also \cite{renaut2020fast}. For square $A$, the theoretical ratio of the cost of SB implemented with $A$ as compared to the use of $\tilde{A}_k$, as given in $(i)$ and $(ii)$ in column $4$ of \Cref{Tab:costs}, is  
\begin{equation}\label{eq:ratio of cost}
    \text{SB Speed up} = \frac{N+2P}{2kn_1}.
\end{equation}  
In practice, we have $P\approx N$  and $k \ll R\leq n_1$. Therefore, \cref{eq:ratio of cost}, which gives the estimated maximum speed up in the SB algorithm with $\tilde{A}_k$ as compared to $A$ is always greater than $1$ which indicates the effectiveness of using the KP approximation as 
compared to the direct  use of $A$, when the augmented system \cref{eq:agum x} is solved using \Cref{alg:CGLS}. Suppose $P=N$ then we get $3N/2kn_1$ but for square $A$ assume we also have $n_1=\sqrt{N}$. Hence we get roughly $3\sqrt{N}/2k$, so the benefit scales roughly as $(2/3)(k/\sqrt{N})$. For larger $N$  and small $k$ we expect dramatic cost reduction with $\tilde{A}_k$ in the implementation of the SB algorithm.  On the other hand, the cost of finding $\tilde{A}_k$ is significant. Under the same assumptions that $A$ is square, column $3$ is negligible as compared to column $2$,  $\mathcal{O}(Nk)$ as compared to $\mathcal{O}(N^2k)$. Therefore, we consider  
\begin{equation}\label{eq:ratio of alg cost}
    \text{Alg Speed up} = \frac{4(N+2P)N\iota_{total}}{4\rho N^2k_p}, 
\end{equation} 
the ratio between column $4$ and column $2$,  where $\rho=1$ or $2$ for EGKB and RSVD, respectively. \Cref{eq:ratio of alg cost}   scales like $ (3/\rho)(\iota_{total} / k_p) \ge 1.5 \iota_{total} / k_p >1$, when $k_p$ is small relative to the total number of iterations of CGLS.  Notice if we need to solve more than one problem, then the  significance of the relatively high cost for finding $\tilde{A}_k$ is reduced, as indicated by \cref{eq:ratio of cost}.  

It is also important to consider the memory demands of the presented algorithms. First of all we note that it is clear that the memory demands for $\tilde{A}_k$ are much less than storing the matrix $A$. While $\tilde{A}_k$  has just $k(m_1n_1+m_2n_2)$ entries,  $A$ has $MN=m_1m_2n_1n_2$ entries. As an example, suppose that $A$ is square, with $M=N=m_1^2=n_1^2$. Then the storage for $\tilde{A}_k$ is just $2kM \ll M^2$ when $k$ is small relative to $M$. But, the use of the KP approximation to $A$ does present a memory challenge if we require to maintain both $\mathcal{R}(A)$ and $A$ in memory concurrently. In particular, $\mathcal{R}(A)$ is needed to find $\tilde{A}_k$, but is not needed elsewhere. Hence, in our experiments we overwrite $A$ in memory by $\mathcal{R}(A)$, and then also delete $\mathcal{R}(A)$ after the SVD in \cref{eq:SVDofR(A)} has been obtained. In our discussion, we have not discussed the costs of reordering of $A$ to find $\mathcal{R}(A)$. We will find that this becomes significant as the size of $A$ increases. 
\subsection{Single Precision to improve the computational cost of \texorpdfstring{$\tilde{A}_k$}{}}\label{seubsec:SP} 
Using SP for large computations is attractive because it reduces  storage and floating point costs as compared to use of double precision. Further, it has been shown that there is an advantage to mixing two precision formats for obtaining accurate results in a reasonable time \cite{higham2022mixed}. In this work, we adopt a mixed precision using both single and double precision. We investigate the use of  SP for the construction of $\tilde{A}_k$ and double precision for the SB reconstruction. The possible consequence of the accumulated errors in SP is that EGKB may suffer from a loss of orthogonality even with one-sided reorthogonalization. To improve the accuracy of the results obtained with EGKB(SP), we can find the KP approximation with full reorthogonalization in place of one-sided reorthogonalization. In relation to \Cref{alg:EGKB}, this means that for SP the columns $\bfs_j$ in step $3$  can be reorthogonalized against all previous columns at each step $j$. Although full reorthogonalization increases the cost of EGKB(SP) algorithm, the dominant cost of the EGKB algorithm, given in \cref{eq:EGKB cost1}, is not impacted. It is important to emphasize that one-sided reorthogonalization is sufficient for EGKB, but the orthogonality of the columns in the EGKB(SP) are maintained using full reorthogonalization.

Using SP reduces the computational cost by a factor of $2$ compared to the use of the standard double precision  \cite{higham2022mixed}. Specifically, the computational costs of EGKB(SP) and RSVD(SP) are significantly lower than the costs  presented in $(i)$ and $(ii)$ in \Cref{Tab:costs}. In addition, forming the term $A_{x_i}$ in the KP approximation is expected to be much faster because $\bar{U}_i$ and $\bar{\sigma}_i$ are stored in SP objects. We use the Matlab function 
$\texttt{single}$ to generate the results in SP.

To conclude, we use mixed precision in this work, where $\tilde{A_k}$ is constructed super fast from SP and the solution is obtained accurately from double precision SB. 
\section{Numerical Experiments}\label{sec:Numerical examples}
Here we present numerical experiments designed to test both the performance of the EGKB(SP) and RSVD(SP) for finding the KP approximation  $\tilde{A}_k$   of $A$, as well as the performance of the $\tilde{A}_k$  for the solution of the TV problem using the SB algorithm. For small-scale problems, EGKB(SP) and RSVD(SP) are contrasted with the use of double precision for the computations in EGKB and RSVD. In addition, we contrast the use of the approximate operators $\tilde{A}_k$ with the use of TSVD $A_c$ for $A$, where the choice of $c$ is described in \cref{subsec:testEGKB}, 
and with $A$  for the direct solves within the SB algorithm. 

We first  present experiments that examine the use of the EGKB(SP) for the determination of the SVD of $\mathcal{R}(A)$ where we may not know in advance how to choose the desired number of terms in the sum approximation. We then describe the choice for the convergence and regularization parameters, and the measures for the performance of the image deblurring algorithms, with the resulting the $\tilde{A}_k$. Results on smaller-scale simulations are used for detailed examination of the parameter choices $\lambda$ and $\beta$ required in the SB algorithm.  These parameter choices are used for the large scale examples, for which it is not practical to search over parameter space.

All the presented numerical results were performed using Matlab Version 2024a on a 2023 Mac Studio with M2 max chip and 64GB of memory. The software will be available on request to the authors.

\subsection{The Experimental Design}\label{subsec:experdesign}
For a square matrix $A$ (i.e., $M=N$), we set $m_1 = m_2 = n_1 = n_2$ to obtain $\mathcal{R}(A)$. We use the modified Gram-Schmidt process to reorthogonalize the columns of $Q_k$ with the EGKB(SP) and EGKB algorithms. Note that reorthogonalizing the columns of $S_{k+1}$ in EGKB(SP) is also required. For the RSVD algorithm we use a single power iteration and chose the number of terms $k$ for the approximation using the estimate obtained from the EGKB(SP) algorithm.

\subsection{Evaluating the KP sum approximation using EGKB(SP) and RSVD(SP)}\label{subsec:testEGKB}
\begin{figure}[ht!]
\centering
\subfloat[Mild blur ]{\label{Fig:Mild Blur}\includegraphics[width=.20\textwidth]{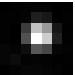}}
\
  \subfloat[Medium blur]{\label{Fig:Medium Blur}\includegraphics[width=.20\textwidth]{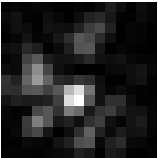}}
\caption{The non separable speckle PSFs used in \Cref{ex:Blur} and obtained from \cite{gazzola2019ir}, with a zoom from the original axes of the PSFs. The sizes of \cref{Fig:Mild Blur,Fig:Medium Blur} are $7 \times 7$ pixels and $15 \times 15$ pixels, respectively.}
\end{figure}

We examine the use of the EGKB(SP) and RSVD(SP) algorithms for $\mathcal{R}(A)$ in constructing $\tilde{A}_k$, $k<R$, as described in \Crefrange{subsec:reorder A}{subsec:RVD}. The EGKB algorithm relies on a starting vector. In our experiments we use  a random vector selected from a normal distribution with mean $0$ and variance $1$. Through our experiments we have found that the choice of starting vector does not impact the rate of convergence of the algorithm.  For the EGKB we also evaluate, for small-scale cases,  the use of the new stopping test given in \cref{eq:stopnu} as a mechanism to 
 estimate a suitable number of terms $k$ for the TSVD approximation. 
We demonstrate the test to find $k$ for the EGKB(SP) and EGKB and the impact on the approximation of the spectrum as compared to the RSVD(SP) and RSVD algorithms using \Cref{ex:Blur}.
\begin{exmp}\label{ex:Blur}
    We consider two blurring operators from \cite{gazzola2019ir}. The image (Mild blur) shown in \Cref{Fig:Mild Blur} is a mild non separable speckle PSF of size $100 \times 100$ with zero boundary conditions. In contrast to the mild blur, the image (Medium blur) shown in \Cref{Fig:Medium Blur} is a medium non separable speckle PSF of size $128 \times 128$ with reflexive boundary conditions. The provided PSFs are spatially invariant blurs, and the parameters that define the blurring matrices $A_{\text{Mild}}$ and $A_{\text{Medium}}$ are given in \Cref{Tab:Blur}. 
\end{exmp}

\begin{figure}[ht!]
 \centering
   \subfloat[Mild Blur: $\nu_j=\sqrt{\zeta_j \zeta_{j-1}}$]{\label{Fig:prodc Mild Blur}\includegraphics[width=.40\textwidth]{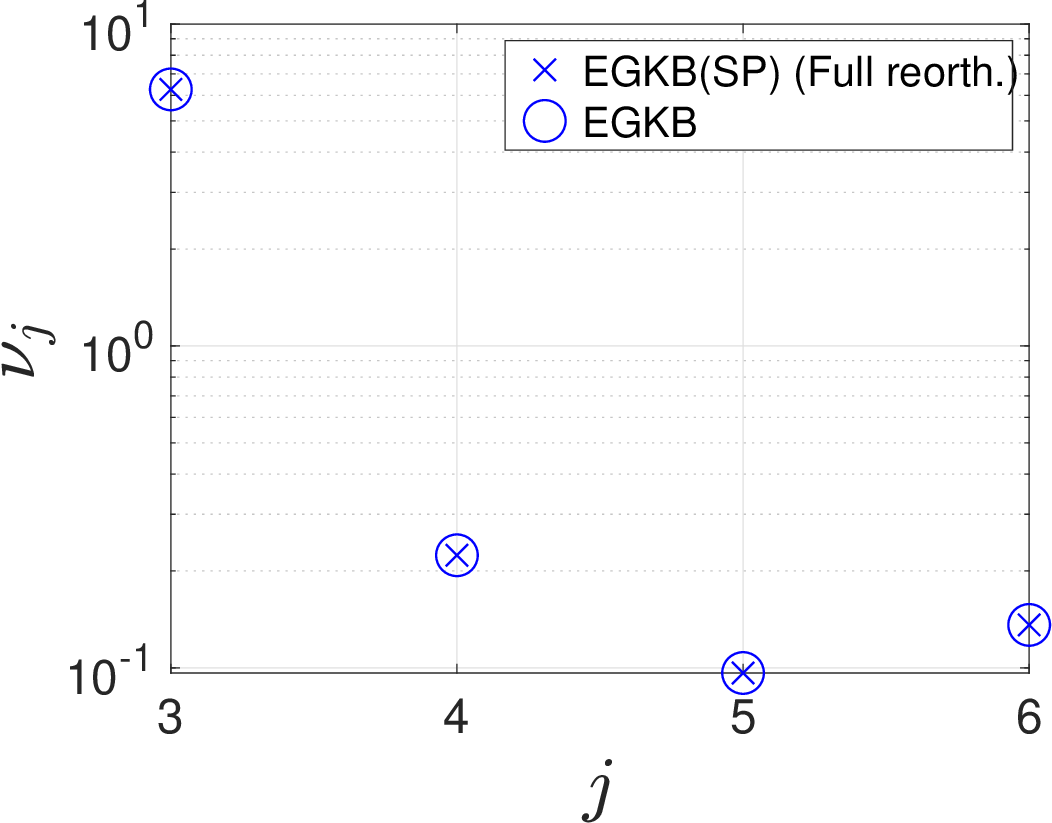}}
    \hfil
      \subfloat[Medium blur: $\nu_j=\sqrt{\zeta_j \zeta_{j-1}}$]{\label{Fig:prodc Medium Blur}\includegraphics[width=.40\textwidth]{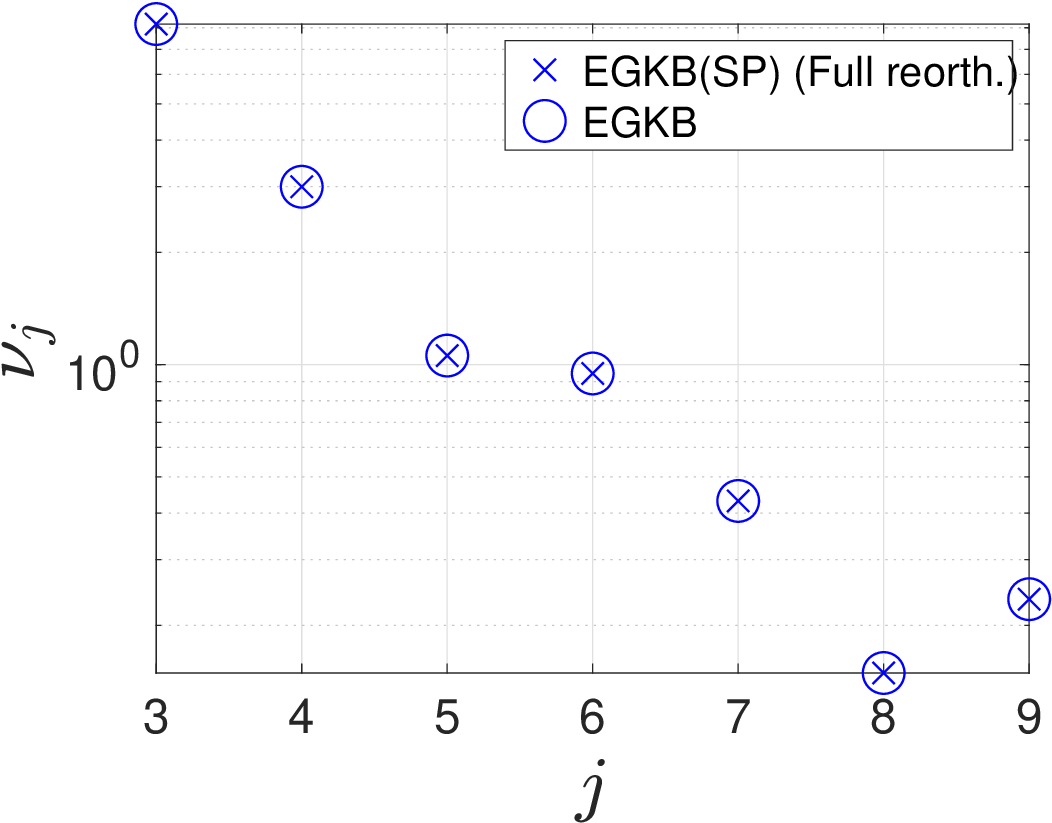}}
     \\
  \subfloat[Mild Blur: SVs]{\label{Fig:TSVD EGKB RSVD of R Mild}\includegraphics[width=.40\textwidth]{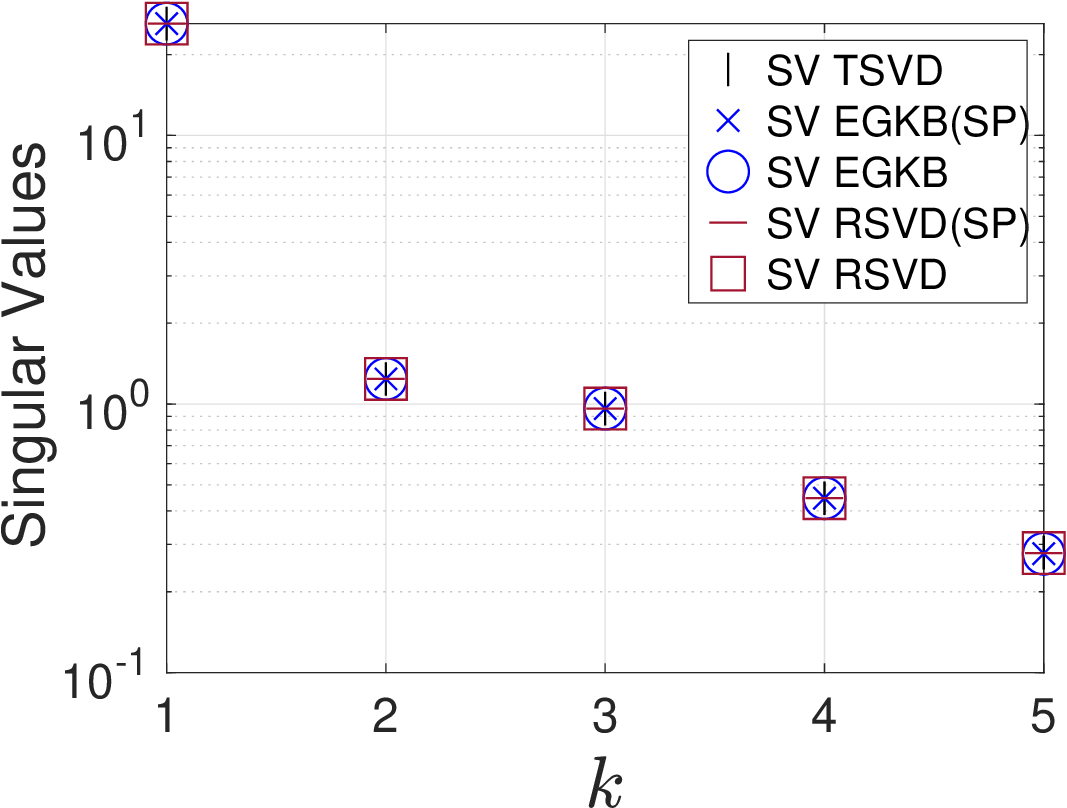}} \hfil
  \subfloat[Medium blur: SVs]{\label{Fig:TSVD EGKB RSVD of R Medium}\includegraphics[width=.40\textwidth]{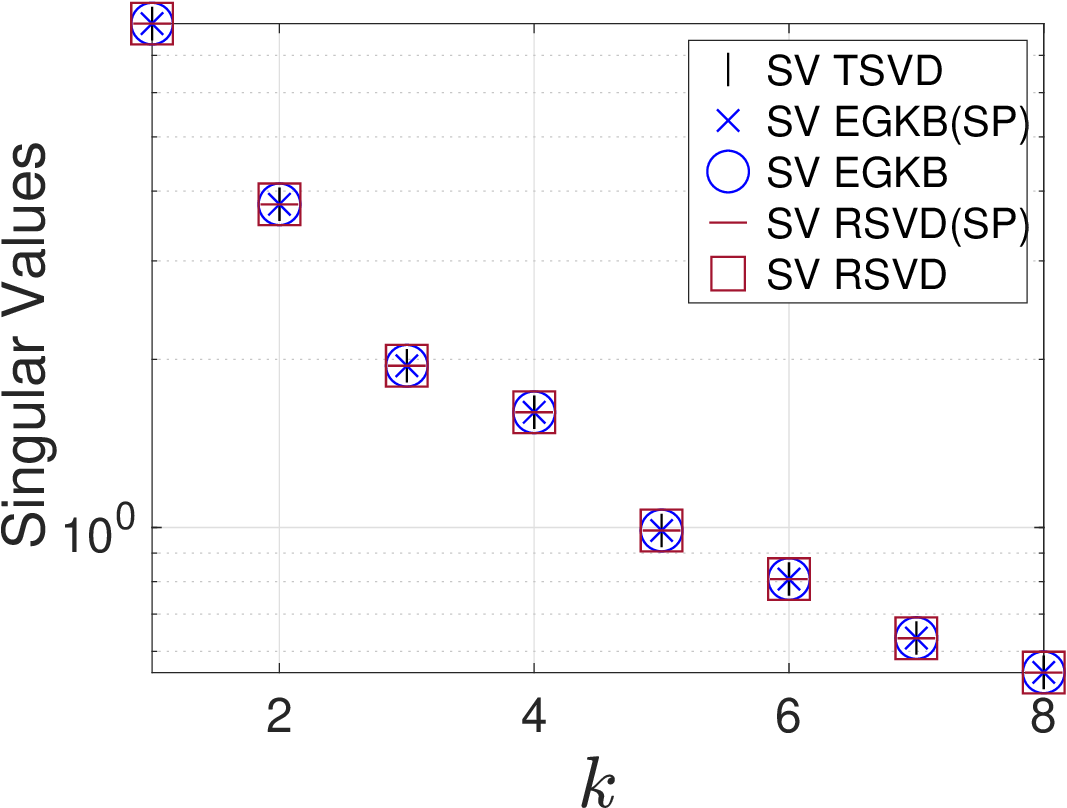}}
 \caption{ \Cref{Fig:prodc Mild Blur} and \Cref{Fig:prodc Medium Blur} show the monotonic decrease of $\nu_j$, as calculated at step $8$ of \Cref{alg:EGKB} using the  geometric means of $\zeta_j$,  until $j=5$ and $j=8$ for \Cref{alg:EGKB}  applied to $\mathcal{R}(\texttt{single}(A_{\text{Mild}}))$ and $\mathcal{R}(\texttt{single}(A_{\text{Medium}}))$, respectively. \Cref{Fig:TSVD EGKB RSVD of R Mild} and \Cref{Fig:TSVD EGKB RSVD of R Medium} show the approximations of the singular values for $\mathcal{R}(A_{\text{Mild}})$ and $\mathcal{R}(A_{\text{Medium}})$ in \Cref{ex:Blur} using the TSVD, EGKB(SP), EGKB, RSVD(SP), and RSVD with the chosen parameters given in \Cref{Tab:Blur}. The black vertical lines, blue crosses, blue circles, dark red horizontal lines, and dark red squares indicate the estimates of the singular values obtained using TSVD, EGKB(SP), EGKB, RSVD(SP), and RSVD, respectively. All presented results with EGKB(SP) are obtained with full reorthogonalization.
 \label{Fig:Blur reults}}
 \end{figure}

\begin{center}
\begin{table}[ht!]
\centering
\caption{\label{Tab:Blur}The parameters that define the blurring matrices in \Cref{ex:Blur}, the number of terms $k$ used for the $\tilde{A}_k$  as determined by the tolerance for $\nu_j$, and the oversampling parameter $p$.}
\begin{tabular}[t]{ c c c c c c c c c c }
\toprule
Blurring Mat &Size &Rank & Condition number& Kronecker rank & $k$ & $p$\\ \midrule
$A_{\text{Mild}}$ &$10000 \times 10000$& $10000$ & $1.5 \times 10^7$& $99$ & $5$& $2$ \\
$A_{\text{Medium}}$ &$16384 \times 16384$& $16384$ &$1.0 \times 10^{8}$ & $127$ & $8$& $4$\\
\bottomrule
\end{tabular}
\end{table}  
\end{center}
\Cref{Fig:prodc Mild Blur,Fig:prodc Medium Blur} show that $\nu_j$, as calculated at step $8$ of \Cref{alg:EGKB} using the  geometric means of $\zeta_j$,  decreases monotonically until $j=5$ and $j=8$ for \Cref{alg:EGKB}  applied to $\mathcal{R}(\texttt{single}(A_{\text{Mild}}))$ and $\mathcal{R}(\texttt{single}(A_{\text{Medium}}))$, respectively. Full reorthogonalization at steps $3$ and $6$ of \Cref{alg:EGKB} is  applied to $\mathcal{R}(\texttt{single}(A_{\text{Mild}}))$ and $\mathcal{R}(\texttt{single}(A_{\text{Medium}}))$. The Krylov space is extended to size $k+p$ with $p_{\text{Mild}}=2$ and  $p_{\text{Medium}}=4$, respectively,  to obtain an estimate of the TSVD for the enlarged space. This estimated SVD is truncated  to $k_{\text{Mild}}=5$ and $k_{\text{Medium}}=8$ terms, respectively, for use in \cref{eq:KPsum}.
In \Cref{Fig:TSVD EGKB RSVD of R Mild,Fig:TSVD EGKB RSVD of R Medium}, we contrast the first $k$ singular values obtained using \Cref{alg:EGKB,alg:RSVD} with the TSVD estimate of these singular values, for $\mathcal{R}(A_{\text{Mild}})$ and $\mathcal{R}(A_{\text{Medium}})$, respectively. \Cref{Fig:TSVD EGKB RSVD of R Mild,Fig:TSVD EGKB RSVD of R Medium} illustrates that good approximations are obtained for the dominant spectrum with the chosen $k$ terms in each case. 

To evaluate the use of the EGKB(SP) and RSVD(SP) algorithms in finding a suitable approximation to $\mathcal{R}(A)$, and hence a good $\tilde{A}_k$ approximating $A$,  we consider the distance of $\bar{\mathcal{R}}_k(A)$ to $\mathcal{R}(A)$, measured in the Frobenius norm, as well as  the error in $\tilde{A}_k$ as an approximation to $A$.  
The baseline is given by the exact $\tilde{A}_k$ using the exact TSVD for $\mathcal{R}(A)$, rather than the approximate $k$-term expansion, that  is closest to $A$  over all other possible Kronecker product representations with $k$ terms, as measured in the Frobenius norm, \cite[Theorem 12.3.1]{golub2013matrix}. We also consider the distance to $A$ for $\tilde{A}_k$ as compared to the best rank $c$ approximation to $A$. In particular, for a general matrix  $A$ with rank $r$, then the TSVD matrix 
$A_c$, with $c<r$, is a best rank $c$ approximation to $A$ if we have $\|A-A_c\|_F=\sqrt{\sum_{i=c+1}^r \sigma_i^2}$,  where $\sigma_i$ are the singular values of $A$. In our comparison $c$ is found such that $\|A-A_c\|_F$ is close to the minimum of the errors in $\tilde{A}_k$ to $A$  obtained using the SP and double precision computations in EGKB and RSVD. This is easy to do by using the square root of the cumulative sums of the vector with entries $\sigma_i^2$, $i=r \dots 1$ and comparing with the desired minimum value.

The errors in $\bar{\mathcal{R}}_k(A)$ to $\mathcal{R}(A)$ and 
in the $\tilde{A}_k$ to $A$, obtained using EGKB(SP), EGKB, RSVD(SP), and RSVD for the parameters $k$ and $p$  given in \Cref{Tab:Blur}, for the matrices $A_{\text{Mild}}$ and $A_{\text{Medium}}$, are summarized in \Cref{Tab:Blur Accur}.
We also present the error for the  TSVD  $A_c$ to $A$, for which $\sqrt{\sum_{i=c+1}^r \sigma_i^2}$ is chosen as the minimum of the errors in the $\tilde{A}_k$ to $A$. We find that the accuracy of the $\tilde{A}_k$ for $A_{\text{Mild}}$ and $A_{\text{Medium}}$ is comparable to $A_{6817}$ and $A_{8072}$, namely $c=6817$ and $8072$, respectively, in each case near $70\%$ of the rank of the matrix. Equivalently, to get comparable accuracy in terms of the approximation to $A$ as measured in the Frobenius norm, we need to take a TSVD with a high number of terms. Finally, the approximation errors for the EGKB(SP), EGKB, RSVD, and RSVD(SP) are comparable in each case, and are independent of the use of SP for the formation to $\tilde{A}_k$.

 \begin{center}
\begin{table}[ht!]
\centering
\caption{\label{Tab:Blur Accur}The relative errors in 
$\bar{\mathcal{R}}_k(A_{\text{Mild}})$ and  $\bar{\mathcal{R}}_k(A_{\text{Medium}})$ to $\mathcal{{R}}(A)$,  and the relative errors in $\tilde{A}_{k_{\text{Mild}}}$ and $\tilde{A}_{k_{\text{Medium}}}$ to $A$ for mild and medium blur, respectively, as given in \Cref{ex:Blur}. All presented results with EGKB(SP) are obtained with full reorthogonalization. (SP) denotes  results obtained with single  precision.} 
\begin{tabular}[t]{c l c c c c c c c c c c}
\toprule
Matrix&Method&$\frac{\|\mathcal{R}(A)-{\bar{\mathcal{R}}}_k(A)\|_F}{\|\mathcal{R}(A)\|_F}$& $\frac{\|A-\tilde{A}_{\text{k}}\|_F}{\|A\|_F}$&$c$&$\frac{\|A-A_{c}\|_F}{\|A\|_F}$\\ \midrule
\multirow{4}{*}{$A_{k_{\text{Mild}}}$}& EGKB(SP)&$0.0115$&$0.0117$&\multirow{4}{*}{$6817$}&\multirow{4}{*}{$0.0116$}  \\
&EGKB&$0.0117$&$0.0117$ \\ 
&RSVD(SP)&$0.0115$&$0.0116$ \\ 
&RSVD&$0.0116$&$0.0116$ \\ \hline
\multirow{4}{*}{$A_{k_{\text{Medium}}}$}& EGKB(SP)&$0.0644$&$0.0666$&\multirow{4}{*}{$8072$}&\multirow{4}{*}{$0.6654$} \\
&EGKB&$0.0666$&$0.0666$ \\ 
&RSVD(SP)&$0.0646$&$0.0667$\\ 
&RSVD&$0.0667$&$0.0667$\\ 
\bottomrule   
\end{tabular}
\end{table}  
\end{center} 
It can be seen from \Cref{Fig:TSVD EGKB RSVD of R Mild,Fig:TSVD EGKB RSVD of R Medium,Tab:Blur} that it is efficient to use the decay of $\nu_j$ to select the number of terms $k$ for $\tilde{A}_k$ found using EGKB(SP). In the rest of this paper, $k$ is determined with EGKB(SP).  For the presented problems, the error in  $\tilde{A}_k$ to $A$ is inherited from the error in $\bar{\mathcal{R}}_k(A)$, this is seen by comparing the results in the third column against the results in the fourth column of \Cref{Tab:Blur Accur}.

\subsubsection{Tests of Convergence for the SB Algorithm}\label{subsubsec:conv}
We note that the convergence of the SB algorithm depends not only on the outer iteration $\ell$ but also the inner iteration $\iota$ in the CGLS, and thus depends on two convergence tolerances $\tau_{SB}$ and $\tau_{CGLS}$, respectively.

\noindent{\bf Convergence of CGLS:} 
For the convergence of the CGLS algorithm in our small scale simulations we chose to look at the relative change in the solution in the CGLS,  as given from the steps 4 and 5 in \Cref{alg:CGLS}, by 
\begin{equation}\label{eq:RC}
    \text{RC}_{CGLS}(\bfx^{(\iota)})=\frac{\|\mu^{(\iota)} \bfw^{(\iota)}\|_2} {\|\bfx^{(\iota)} \|_2},
\end{equation}
corresponding to the relative change, ($\text{RC}_{CGLS}$), in the solution at iteration $\iota$. We observe that an acceptable approximation along the search direction $\bfw^{(\iota)}$ can be obtained when the  $\text{RC}_{CGLS}$, as given in \cref{eq:RC}, is below the tolerance $\tau_{CGLS}=10^{-4}$. 

\noindent{\bf Convergence of the SB:}
To determine the convergence of the SB algorithm, we can look at the convergence of the relative change in the solutions over two consecutive iterations, specifically we record 
\begin{equation}\label{eq:RCSB}
        \text{RC}_{SB}(\bfx^{(\ell)})=\frac{\|\bfx^{(\ell)}- \bfx^{(\ell-1)}\|_2} {\|\bfx^{(\ell-1)}\|_2}, 
\end{equation}
at each iteration. We consider the solution to be converged when $\text{RC}_{SB}(\bfx^{(\ell)}) < \tau_{SB}$. The choice of convergence tolerance is significant to the quality of the solution, and  should be chosen based on the noise in the problem. Specifically, when the noise in the data is higher, we stop the iteration with a small tolerance compared to the case with lower noise. This test is frequently used for the examination of convergence when no true data is available. Here we have true data but we also use this test to validate the simulation.

\noindent{\bf Maximum Number of Iterations:}
Stopping criteria are imposed for all iterations by introducing a limit on the number of iterations of the algorithm. 
For both the outer iterations with the SB algorithm and the inner iterations with the CGLS algorithm a limit on the number of iterations is  imposed via $L_{max}$, and $I_{max}$, respectively. The number of SB iterations should be relatively small as compared to the size of the problem being solved.  
Consequently, we chose a quite generous upper bound $L_{max}=50$.  The CGLS algorithm solves, in contrast, a large system of equations for both isotropic and anisotropic regularizers, and it is therefore reasonable to pick $I_{max}$ large enough that the iteration captures the dominant spectrum of the underlying problem. Here we chose $I_{max}=100$.

\subsection{The Regularization Operators}\label{subsec:regmat}
For our experiments, the anisotropic and isotropic regularization matrices $L_x$ and $L_y$ in \cref{eq:x in 2d}, are given by $L$ defined by 
\begin{equation}\label{eq:deriv}
L=  \frac{1}{2} \begin{bmatrix}
       1&-1 &&&&0\\
       &1&-1&&&&\\
       &&1&-1&&\\
       &&&\ddots&\ddots&\\
       0&&&&1&-1
   \end{bmatrix}\in \mathbb{R}^{m_2 \times n_2}, 
\end{equation}
for the appropriate size in each of the $x$ and $y$ directions. 
We assume that the images are square, thus $L_x=L_y$ and for simplicity we set $\lambda_x=\lambda_y = \lambda$ and $\beta_x=\beta_y=\beta$ in \cref{eq:agu A}. Where the blur operator is known to be spatially variant in the $x$ and $y$ directions, it would be more suitable to use different regularization parameters in each direction.  

\subsubsection{Regularization parameters}
For the selections of the regularization parameters $\lambda$ and $\beta$ in our simulations, we estimate suitable parameters using the ratios $\gamma$ which define the shrinkage and $\lambda$ relative to the noise level. For $\lambda$, we sweep over $100$ logarithmically spaced points in the $(\lambda, \beta)$ domain, where $\beta$  is implicitly defined as  $\beta=\lambda^2/\gamma$, for a fixed $\gamma$. The near optimal parameters are fixed values for each reconstruction and are only found to illustrate the effectiveness of the methods. For practical problems estimating  suitable regularization parameters remains a challenge but is guided by assumed knowledge of the noise in the data, a reasonable requirement for measured data. Examples of approaches to optimally find $\beta$ and/or $\lambda$ are discussed and reviewed in \cite{CHUNG2023297,sweeney2024parameterselectiongcvchi2}.
\subsubsection{Performance Metrics}\label{subsubsec:perfmeas}
The Signal-To-Noise ratio (SNR) is used to quantify the noise in the data, usually expressed in decibels (dB), and is given by
\begin{align}\label{eq:SNR}
    \text{SNR}(\bfbtrue, \bfb)= 10 \ \text{log}_{10} \ \frac{\| \bfbtrue\|_2^2}{ \|\bfbtrue - \bfb\|_2^2}= 10 \ \text{log}_{10} \ \frac{\| \bfbtrue\|_2^2}{ \|\bfeta \|_2^2}.
\end{align} 
Noisy data $\bfb$ is obtained by adding a Gaussian white noise vector to $\bfbtrue$ according to $ \bfb= \bfbtrue + \bfeta$, where $\bfeta =\eta (\|\bfbtrue\|_2/\|\bfzeta\|_2) \bfzeta$, 
for noise level $\eta$. Here $\bfzeta=\texttt{randn}(M,1)$ where \texttt{randn} is the Matlab function that generates a vector of normally distributed random elements with mean zero and variance of one.

The relative error in the resulting image is recorded as 
\begin{equation}\label{eq:RE}
    \text{RE}(\bfx^{(\ell)})=\frac{\|\bfx^{(\ell)}- \bfxtrue\|_2} {\|\bfxtrue\|_2},
\end{equation}
where we assume $\bfxtrue$ is available for these simulations, and we also record the improved signal to noise ratio (\text{ISNR}) 
\begin{align}\label{eq:ISNR}
    \text{ISNR}(\bfx^{(\ell)}) = 20\log_{10}\left(\frac{\|\bfb-\bfxtrue\|_2}{\|\bfx^{(\ell)}-\bfxtrue\|_2}\right).
\end{align}
As $\bfx^{(\ell)}$ becomes closer to $\bfxtrue$ the \text{ISNR} should improve, whereas the \text{RE} should decrease.

For each simulation we record the \text{RE}, \text{ISNR}, and the \text{RC} in the solution in the  SB, with increasing outer iteration $\ell$.  
In the simulations to determine the optimal choices for the convergence parameters we do not record the timings, these are only recorded after we have determined the choices of all parameters. For the reconstruction obtained with the $\tilde{A}_k$, we record the timing for each step needed for constructing $\tilde{A}_k$ separately from the timing of the SB inversion. Specifically, we record the timing to find the reordering $A$ to $\mathcal{R}(A)$, the TSVD approximation $\bar{\mathcal{R}}_k(A)$ to $\mathcal{R}(A)$, and the terms in the KP approximation $\tilde{A}_k$ to $A$. We denote the timing to find the terms in the KP approximation by Terms($\tilde{A}_k$).

\begin{figure}[ht!]
 \centering
  \subfloat[True image]{\label{Fig:True pattern2}\includegraphics[width=.20\textwidth]{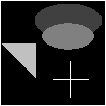}}
   \subfloat[Contaminated]{\label{Fig:blurred Patteren2}\includegraphics[width=.20\textwidth]{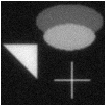}}
 \subfloat[ANI $A$ \label{Fig:Dir Ani}]{\includegraphics[width=.20\textwidth]{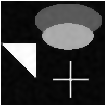}}
  \subfloat[ISO $A$ \label{Fig:Dir iso}]{\includegraphics[width=.20\textwidth]{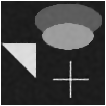}}  
  \\
   \subfloat[ANI EGKB(SP) \label{Fig:EGKB(SP) Ani}]{\includegraphics[width=.20\textwidth]{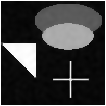}} 
      \subfloat[ANI EGKB \label{Fig:EGKB(DP) Ani}]{\includegraphics[width=.20\textwidth]{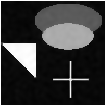}} 
    \subfloat[ANI RSVD(SP)  \label{Fig:RSVD Ani SP}]{\includegraphics[width=.20\textwidth]{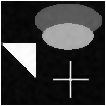}} 
       \subfloat[ANI RSVD  \label{Fig:RSVD Ani DP}]{\includegraphics[width=.20\textwidth]{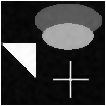}} 
    \\
 \subfloat[ISO EGKB(SP) \label{Fig:EGKB iso SP}]{\includegraphics[width=.20\textwidth]{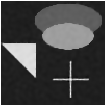}} 
 \subfloat[ISO EGKB \label{Fig:EGKB iso DP}]{\includegraphics[width=.20\textwidth]{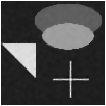}} 
 \subfloat[ISO RSVD(SP) \label{Fig:RSVD iso SP}]{\includegraphics[width=.20\textwidth]{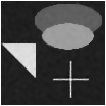}} 
  \subfloat[ISO RSVD \label{Fig:RSVD iso DP}]{\includegraphics[width=.20\textwidth]{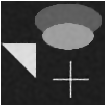}} 
 \caption{A true image of $100 \times 100$ pixels \Cref{Fig:True pattern2} is blurred using a mild non separable speckle PSF of size  $100 \times 100$ \Cref{Fig:Mild Blur} from \cite{gazzola2019ir} and $7\%$ Gaussian noise, which has an approximate SNR$ =23$ for the problem in \cref{ex:Pattern2}. SB solutions obtained with $A$ and the $\tilde{A}_k$ at convergence for \Cref{Fig:blurred Patteren2} using $\lambda=0.1150$ and $\beta=0.0066$. Here the approximations to $\tilde{A}_k$, denoted with (SP), are obtained  for $R(\texttt{single}(A))$.}\label{Fig:pattern2} 
 \end{figure}

\noindent{\bf Acronyms for the simulation:} We use the acronyms ANI and ISO to denote the anisotropic and isotropic SB algorithm. These acronyms are coupled with $A$, EGKB(SP), EGKB, RSVD(SP), RSVD to indicate the application of the SB algorithm with $A$ directly and $\tilde{A}_k$ obtained from EGKB(SP), EGKB, RSVD(SP), and RSVD, respectively. 
 \begin{exmp}\label{ex:Pattern2}
     We consider the \texttt{pattern2} image of size $100 \times 100$ from \cite{gazzola2019ir}, shown in \Cref{Fig:True pattern2}. We blur this image with a mild non separable speckle PSF given in \Cref{Fig:Mild Blur}. The provided PSF is a spatially invariant blur. Then we use $7 \%$ white Gaussian noise, corresponding to signal to noise ratio (SNR) of approximately $23$, to display the noise contaminated \texttt{pattern2} image in \Cref{Fig:blurred Patteren2}.
       The parameters that define the blurring operator matrix and the number of terms $k$ used for the $\tilde{A}_k$ with the oversampling parameters $p$, as determined by the tolerance for $\nu_j $, are given in the first row of \Cref{Tab:Blur} in \cref{subsec:testEGKB}. 
 \end{exmp} 
   For the parameters in the SB algorithm, we use a convergence tolerance $\tau_{SB} =10^{-3}$.  In addition, we set $\gamma=2$ and swept $\lambda$ across $[10^{-2}, 10]$ to obtain the near optimal parameters. We note that it is efficient to use the parameters resulting from ANI EGKB(SP) for all cases. This observation is justified for the following reasons: (i) the errors in the $\tilde{A}_k$ to $A$ obtained using the EGKB(SP), EGKB,  RSVD(SP), and RSVD are comparable, as given in \Cref{Tab:Blur Accur}, and they both yield almost indistinguishable approximations to the first $k$ singular values of $R(A)$, as shown in \cref{Fig:TSVD EGKB RSVD of R Mild}. Thus, we expect that the approximation  of $\tilde{A}_k$ to $A$ with either EGKB(SP) or RSVD(SP) provide  $\tilde{A}_{k}$ with comparable spectral properties, which influences the resulting choices for $\lambda$.  (ii) Because the error in the approximation obtained with EGKB(SP) is relatively small, as shown in \Cref{Tab:Blur Accur}, the parameters used  in conjunction with $\tilde{A}_k$ can also be reasonably used for the reconstruction with matrix $A$. (iii) We see in \Cref{sec:TVReg} that the only difference between the ANI SB and the ISO SB algorithms is the  $\ell_1$ term, impacting the calculation of $\bfd_x$ and $\bfd_y$. We anticipated that for the $\ell_1$ regularized problem, which enforces sparsity in the solution, and with simple images similar to those provided in this paper, the impact of the difference between the ANI SB and ISO SB algorithms is reduced. Therefore, for the ISO SB we adopt the parameters estimated for use in  ANI SB.  This does not  significantly affect the results, as  has been confirmed in our simulations.  Sweeping over the range of $\lambda$ and $\gamma$, we determine that it is optimal to use $\lambda=0.1150$ and $\beta=0.0066$ (where $\beta$ is derived from $\gamma$) as indicated by these parameters providing a monotonically-decreasing $\text{RC}_{SB}$ given by \cref{eq:RCSB}. 
  The reconstructed \texttt{pattern2} image using the $\tilde{A}_k$  
  is illustrated in \Cref{Fig:pattern2}, and the quantitative results are summarized in \Cref{Tab:parameters Pattern2}.

\begin{figure}[ht!]
 \centering
   \subfloat[ANI  \label{Fig: Ani RE RC SB}]{\includegraphics[width=.45\textwidth]{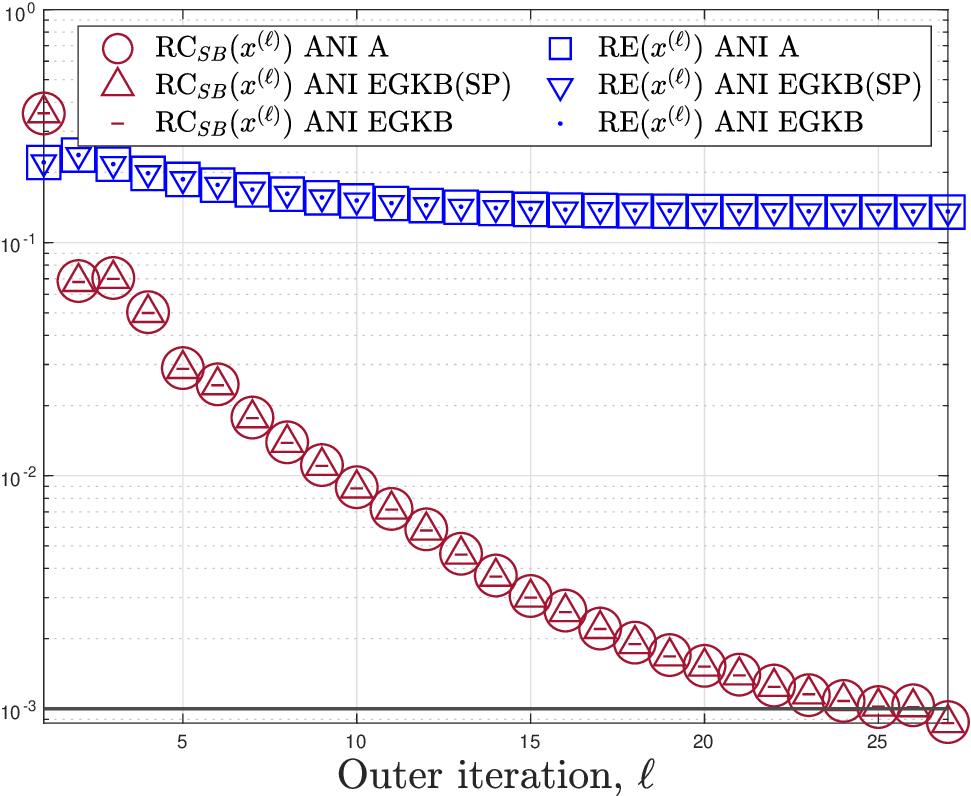}}
   \quad
    \subfloat[ISO \label{Fig: ISO RC CGLS}]{\includegraphics[width=.45\textwidth]{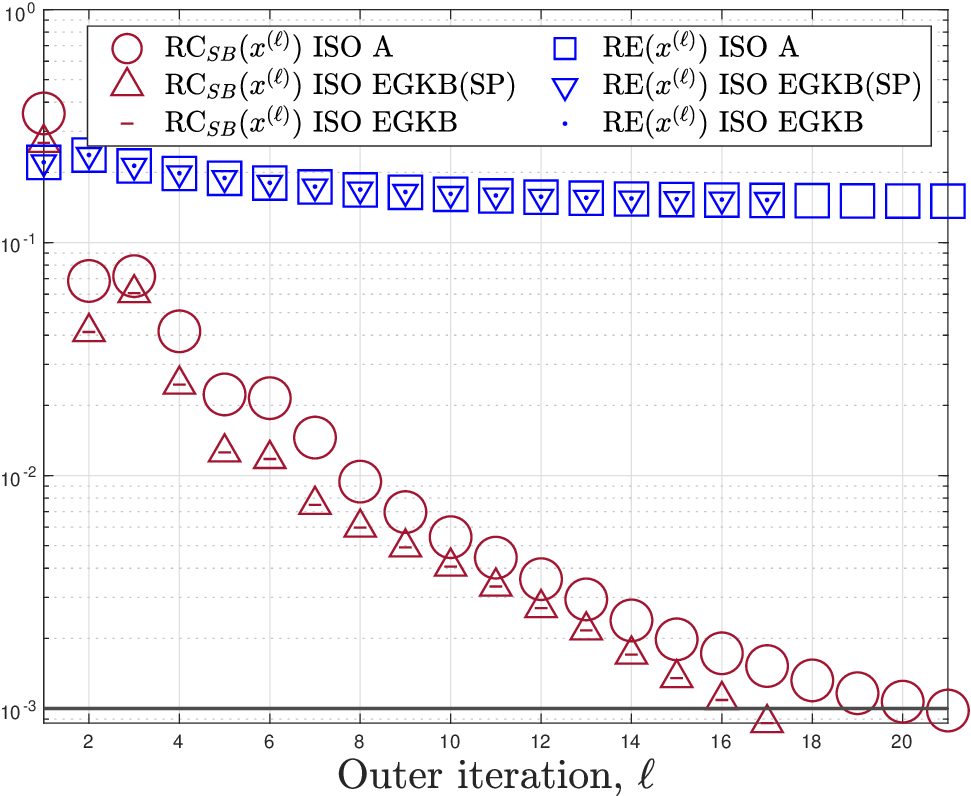}}
 \caption{The corresponding $\text{RC}_{SB}$ and \text{RE} in the outer iterations $\ell$ in the SB algorithm with $A$ directly and with $\tilde{A}_k$ estimated by EGKB(SP) and EGKB for the problem in \cref{ex:Pattern2}. The matrix $A$ is of size $10000\times 10000$. The horizontal line in each plot indicates the convergence tolerance $\tau_{SB}$.}\label{Fig:RE RC pattern2} 
 \end{figure}

 Due to the similarity in most reconstructions, we present the corresponding $\text{RC}_{SB}$ and \text{RE} in the outer iterations $\ell$ in the SB algorithm using ANI $A$, ANI EGKB(SP), and ANI EGKB in \Cref{Fig: Ani RE RC SB}, and  ISO $A$, ISO EGKB(SP), and ISO EGKB in \Cref{Fig: ISO RC CGLS}. The plots confirm that the convergence behavior of $\text{RC}_{SB}$ and \text{RE} is almost independent of the method applied, although the $\text{RC}_{SB}$ for the algorithm using $A$ with the ISO SB converges more slowly than using $\tilde{A}_k$. 

 \begin{table}[ht!]
 \footnotesize
 \caption{\label{Tab:parameters Pattern2} The corresponding \text{RE} and \text{ISNR}, \cref{eq:RE,eq:ISNR}, respectively,  by the SB with $A$ and $\tilde{A}_k$ obtained with EGKB(SP), EGKB, RSVD(SP), and RSVD applied to the problem with $\text{SNR}\approx 23$ in \Cref{ex:Pattern2}. The matrix $A$ is of size $10000\times 10000$. All results are obtained with the SB parameters   $\lambda=0.1150$ and $\beta=0.0066$.} 
 \begin{center}
 \begin{tabular}[t]{l c c c c c c c c c c } 
 \toprule 
 &&&&&\multicolumn{5}{c}{CPU Time}\\ \cline{6-10}  \\[-2ex]
Method&$\ell_{\text{end}}$&$\iota_{total}$&$\text{RE}(\bfx)$&$ \text{ISNR}(\bfx)$ 
 &$\mathcal{{R}}(A)$&$\bar{\mathcal{R}}_k(A)$&Terms($\tilde{A}_k$)& SB&Total  \\ \midrule
ANI $A$      &$27$  &$539$&$0.14$ &$9.27$&$-$&$-$&$-$&$10.62$&$10.62$\\ 
ANI EGKB(SP) &$27$&$464$&$0.14$&$9.22$&$0.37$&$0.14$&$0.0117$&$0.81$&$1.33$\\ 
 ANI EGKB    &$27$&$464$&$0.14$&$9.22$&$0.37$&$0.15$&$0.0117$&$0.75$&$1.28$\\ 
 ANI RSVD(SP)&$27$&$464$&$0.14$&$9.22$&$0.37$&$0.22$&$0.0116$&$0.74$&$1.34$\\ 
 ANI RSVD    &$27$&$464$&$0.14$&$9.22$&$0.37$&$0.35$&$0.0116$&$0.77$&$1.50$\\ 
 ISO $A$     &$21$&$405$&$0.15$&$8.34$&$-$&$-$&$-$&$7.85$&$7.85$\\ 
 ISO EGKB(SP)&$17$&$294$ &$0.15$&$8.23$&$0.37$&$0.14$&$0.0117$&$0.50$&$1.02$\\ 
 ISO EGKB    &$17$&$294$ &$0.15$&$8.23$&$0.37$&$0.15$&$0.0117$&$0.48$&$1.01$\\ 
 ISO RSVD(SP)&$17$&$294$ &$0.15$&$8.23$&$0.37$&$0.22$&$0.0116$&$0.50$&$1.10$\\ 
 ISO RSVD    &$17$&$294$ &$0.15$&$8.23$&$0.37$&$0.35$&$0.0116$&$0.49$&$1.22$\\ 
\bottomrule
 \end{tabular}    
   \end{center} 
 \end{table}

  We see immediately that coupling $\tilde{A}_k$ with the SB algorithm can provide excellent approximations of the solution, similar to that obtained using  $A$,  and that ISO uses fewer outer SB steps. Empirically, the near optimal parameters obtained for the ANI EGKB(SP) are used successfully to approximate all other cases. As reported in column $6$ of \Cref{Tab:parameters Pattern2}  finding $\mathcal{R}(A)$ dominates the cost of constructing $\tilde{A}_k$, and this rearrangement increases with the increasing size of $A$. However, as discussed in \cref{subsec:CostSB}, the computationally most time-consuming  part  in obtaining $\tilde{A}_k$ is to calculate the TSVD approximation to $\mathcal{R}(A)$. Cost estimates indicate that the calculation using the RSVD should be roughly double that of the EGKB. This is confirmed by the results reported in column $7$  of \Cref{Tab:parameters Pattern2}. In this particular problem, the number of terms $k=5$ for $\tilde{A}_k$ and $p=2$ are relatively small compared to the Kronecker rank  of $A$ at $R=99$. Thus, the difference between EGKB(SP) and EGKB is not significant.   Overall, the timing results of the SB algorithm do not indicate a significant difference in cost between using EGKB(SP), EGKB, RSVD,  and RSVD(SP)  for these small cases. 
  
  To compare the theoretical and actual costs of using SB with $A$ or with  $\tilde{A}_k$, we want to apply \cref{eq:ratio of cost} for \Cref{ex:Pattern2}. The relevant parameters given  in \Cref{Tab:parameters Pattern2} are $M=N=P=10000$, with $m_1=n_1=100$, and  $T=N+2P=30000$.  We also  need to assume that  $\iota_{total}$ is roughly the same for both cases, although we can see from \cref{Tab:parameters Pattern2} that $\iota_{total}$ is higher for the direct implementation with $A$. Substituting these parameters in \cref{eq:ratio of cost} predicts that the  cost of SB with $A$ would be roughly $30$ times more expensive than the cost with $\tilde{A}_k$. Practically, using the costs given in column $9$  of \Cref{Tab:parameters Pattern2} we see that  SB with  $A$ is approximately $14$ times, respectively $16$ times, more expensive than $\tilde{A}_k$, when implemented for ANI and ISO regularization, respectively, while the total cost is increased by approximately a factor of approximately $7$ in each case.  Moreover, when considering \cref{eq:ratio of alg cost}, it is clear from these results that in each case $\iota_{total}\gg k$, and the cost of forming $\tilde{A}_k$ is negligible as compared to the cost of using $A$ in the SB algorithm. Assume we want to apply \cref{eq:ratio of alg cost} for ANI EGKB(SP) in \Cref{ex:Pattern2}. Using the relevant parameters given in \Cref{Tab:Blur},   the estimate \cref{eq:ratio of alg cost} shows that forming $\tilde{A}_k$ is $1.5(\iota_{total}/k_p)\approx 99$ times faster than using $A$ in the SB. In practice, the ratio between the costs given through column $9$ and the sum of columns $6-8$ in \Cref{Tab:parameters Pattern2}  shows that $\tilde{A}_k$ is only approximately $20$ times more efficient than using $A$ in the ANI SB. This is still a significant speed-up if it is maintained for large scale problems. 

  While \cref{eq:ratio of cost,eq:ratio of alg cost} overestimate the benefit of the use of $\tilde{A}_k$, we note that they were obtained with a number of assumptions leading to dropping lower order terms, and without considering memory access and usage. These factors, particularly memory access and usage, may become more relevant for larger problem sizes, but obtaining comparative data is  not feasible. The results confirm emphatically that there is significant benefit to usage of $\tilde{A}_k$ in the SB algorithm, which is of particular relevance when the regularized problem is solved repeatedly for multiple right hand sides. Further, as shown in \cref{Fig:RE RC pattern2}, the decrease of the RC and RE with iteration $\ell$ and $\iota$ are comparable across methods and  the obtained results are all comparable in terms of the relative error.  Thus, it should be suitable to use $\tilde{A}_k$ for large-scale problems when using $A$ is infeasible. As an aside,  while the ANI SB and ISO  SB results are close in terms of the achieved RE, the ANI SB results do indicate an improved ISNR as compared to ISO SB, although the difference is imperceptible in \Cref{Fig:pattern2}.

\begin{figure}[ht!]
 \centering
  \subfloat[True image]{\label{Fig:True satellite}\includegraphics[width=.20\textwidth]{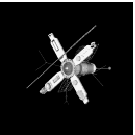}}
   \subfloat[ANI EGKB(SP) \label{Fig:EGKB Ani satellite SP}]{\includegraphics[width=.20\textwidth]{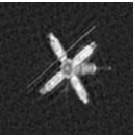}} 
   \subfloat[ANI EGKB \label{Fig:EGKB Ani satellite DP}]{\includegraphics[width=.20\textwidth]{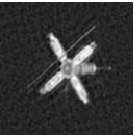}} 
    \subfloat[ANI RSVD(SP)  \label{Fig:RSVD Ani satellite SP}]{\includegraphics[width=.20\textwidth]{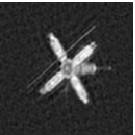}} 
        \subfloat[ANI RSVD  \label{Fig:RSVD Ani satellite DP}]{\includegraphics[width=.20\textwidth]{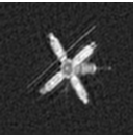}} 
    \\
     \subfloat[Contaminated]{\label{Fig:blurred satellite}\includegraphics[width=.20\textwidth]{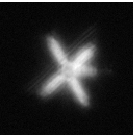}}
 \subfloat[ISO EGKB(SP) \label{Fig:EGKB iso satellite SP}]{\includegraphics[width=.20\textwidth]{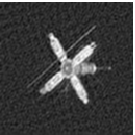}} 
  \subfloat[ISO EGKB \label{Fig:EGKB iso satellite DP}]{\includegraphics[width=.20\textwidth]{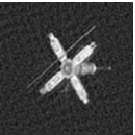}} 
 \subfloat[ISO RSVD(SP) \label{Fig:RSVD iso satellite SP}]{\includegraphics[width=.20\textwidth]{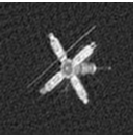}} 
  \subfloat[ISO RSVD \label{Fig:RSVD iso satellite DP}]{\includegraphics[width=.20\textwidth]{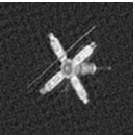}} 
 \caption{A true image of $128 \times 128$ pixels \Cref{Fig:True satellite} is blurred using a medium non separable speckle PSF of size  $128 \times 128$ \Cref{Fig:Medium Blur} from \cite{gazzola2019ir} and $4\%$ Gaussian noise, which has an approximate SNR$ =28$ for the problem in \cref{ex:Pattern2}. SB solutions obtained with $A$ and the $\tilde{A}_k$ at convergence for \Cref{Fig:blurred satellite} using $\lambda=0.0110$ and $\beta=0.0024$.}\label{Fig:satellite} 
 \end{figure}

 \begin{exmp}\label{ex:satellite}
     We consider the \texttt{satellite} image of size $128 \times 128$ from \cite{gazzola2019ir}, shown in \Cref{Fig:True satellite}. We blur this image with a medium non separable speckle PSF given in \Cref{Fig:Medium Blur}. The provided PSF is a spatially invariant blur. Then we use $4 \%$ white Gaussian noise, corresponding to signal to noise ratio (SNR) of approximately $28$, to display the noise contaminated \texttt{satellite} image in \Cref{Fig:blurred satellite}.
 The parameters that define the blurring operator matrix and the number of terms $k$ used for the $\tilde{A}_k$ with the oversampling parameters $p$, as determined by the tolerance for $\nu_j $, are given in the second row of \Cref{Tab:Blur} in \cref{subsec:testEGKB}. 
  \end{exmp}

Due to the error in  this problem as compared to  the error in \cref{ex:Pattern2}, we set $\tau_{SB}=10^{-2}$. Since the PSF significantly affects the quality of the image, we set a smaller shrinkage parameter, $\gamma=0.05$, to maximize the impact of the update \cref{eq:d thresholding ani,eq:d thresholding iso}. We swept $\lambda$ across $[10^{-2},  1]$ to obtain the near optimal parameters. We see that $\lambda=0.0110$ and $\beta=0.0024$ obtained from the ANI EGKB(SP) problem with full reorthogonalization serve as good parameters for all selected methods. The associated computational demands for $A$ in this example make it unattractive. The reconstructed \texttt{satellite} image using the $\tilde{A}_k$  
is illustrated in \Cref{Fig:satellite}, and the quantitative results are summarized in \Cref{Tab:parameters satellite}.
 
 \begin{table}[ht!]
 \footnotesize 
 \caption{\label{Tab:parameters satellite} The corresponding \text{RE} and \text{ISNR}, \cref{eq:RE,eq:ISNR}, respectively,  by the SB with $\tilde{A}_k$ obtained with EGKB(SP), EGKB, RSVD(SP), and RSVD applied to the problem with $\text{SNR}\approx 28$ in \Cref{ex:satellite}. The matrix $A$ is of size $16384\times 16384$. All results are obtained with the SB parameters $\lambda=0.0110$ and $\beta=0.0024$.}
 \begin{center}     
 \begin{tabular}[t]{l c c c c c c c c c c } 
 \toprule 
 &&&&&\multicolumn{5}{c}{CPU Time}\\ \cline{6-10}  \\[-2ex]
Method&$\ell_{\text{end}}$&$\iota_{total}$&$\text{RE}(\bfx)$&$ \text{ISNR}(\bfx)$ 
 &$\mathcal{{R}}(A)$&$\bar{\mathcal{R}}_k(A)$&Terms($\tilde{A}_k$)& SB&Total  \\ \midrule
ANI EGKB(SP) &$13$&$1258$&$0.24$&$7.63$&$2.04$&$0.51$&$0.0009$&$3.50$&$6.05$\\ 
 ANI EGKB    &$12$&$1183$&$0.24$&$7.67$&$2.04$&$0.61$&$0.0012$&$3.29$&$5.94$\\ 
 ANI RSVD(SP)&$13$&$1248$&$0.24$&$7.63$&$2.04$&$0.72$&$0.0013$&$3.52$&$6.28$\\ 
 ANI RSVD    &$12$&$1185$&$0.24$&$7.68$&$2.04$&$0.89$&$0.0014$&$3.31$&$6.24$\\ 
 ISO EGKB(SP)&$5$&$494$&$0.24$&$7.53$&$2.04$&$0.51$&$0.0009$&$1.41$&$3.96$\\ 
 ISO EGKB    &$5$&$494$&$0.24$&$7.53$&$2.04$&$0.61$&$0.0012$&$1.41$&$4.06$\\ 
 ISO RSVD(SP)&$5$&$494$&$0.24$&$7.52$&$2.04$&$0.72$&$0.0013$&$1.41$&$4.17$\\ 
 ISO RSVD    &$5$&$494$&$0.24$&$7.52$&$2.04$&$0.89$&$0.0014$&$1.41$&$4.34$\\ 
\bottomrule 
 \end{tabular}       
   \end{center} 
 \end{table} 
 
From \Cref{Fig:satellite}, we see that the SB coupled with a well defined approximation $\tilde{A}_k$ of $A$ can improve the recovery of the blurred image in \Cref{Fig:blurred satellite}. As recorded in \Cref{Tab:parameters satellite},  the ISO SB converges approximately $7$ iterations earlier than the ANI SB and with similar relative errors, but again the ISNR achieved by the ANI SB is slightly higher. The reported results in column $7$ of \Cref{Tab:parameters satellite} confirm the effectiveness of using SP when forming $\bar{\mathcal{R}}(A)$ compared to the use of double precision. Although the unavoidable loss of accuracy in $\tilde{A_k}$ obtained from EGKB(SP) and RSVD(SP) leads to slightly slower convergence for ANI EGKB(SP) and ANI RSVD(SP), $\ell_{\text{end}}$ is larger, the relative errors in ANI EGKB(SP) and ANI RSVD(SP) are similar to those achieved in double precision. Even though the size of this problem and the number of terms $k$ for $\tilde{A}_k$ are larger than those in \Cref{ex:Pattern2}, the gain in the computation time for the solution required by the SB, as given in column $9$ of \Cref{Tab:parameters satellite}, is reasonable.  Overall, the inclusion of SP for estimating $\tilde{A}_k$ provides solutions with comparable accuracy for both ISO and ANI SB implementations, but demonstrates that we cannot assume that the use of SP will always provide a gain in computational efficiency, particularly if a larger $\ell_{\text{end}}$ is required.

\begin{figure}[ht!]
 \centering
  \subfloat[True image]{\label{Fig:True hst}\includegraphics[width=.20\textwidth]{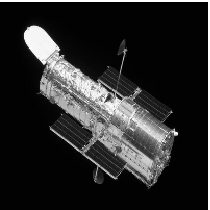}}
 \subfloat[PSF]{\label{Fig:PSF hst}\includegraphics[width=.20\textwidth]{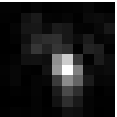}}
  \subfloat[Contaminated]{\label{Fig:blurred hst}\includegraphics[width=.20\textwidth]{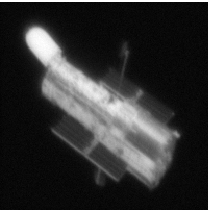}}
  \\
   \subfloat[ANI EGKB(SP) \label{Fig:EGKB Ani hst SP}]{\includegraphics[width=.20\textwidth]{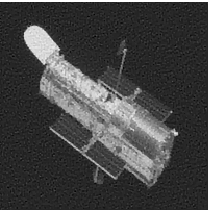}} 
    \subfloat[ANI RSVD(SP)  \label{Fig:RSVD Ani hst SP}]{\includegraphics[width=.20\textwidth]{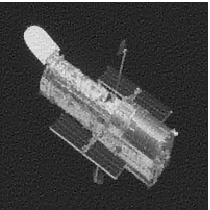}} 
 \subfloat[ISO EGKB(SP) \label{Fig:EGKB iso hst SP}]{\includegraphics[width=.20\textwidth]{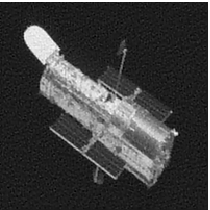}} 
 \subfloat[ISO RSVD(SP) \label{Fig:RSVD iso hstSP}]{\includegraphics[width=.20\textwidth]{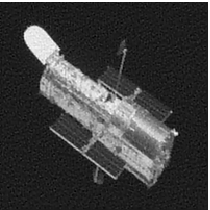}} 
 \caption{A true image of $200 \times 200$ pixels \Cref{Fig:True hst} is blurred using a mild non separable speckle PSF of size  $200 \times 200$ \Cref{Fig:PSF hst} from \cite{gazzola2019ir} and $3\%$ Gaussian noise, which has an approximate SNR$ =30$ for the problem in \cref{ex:hst}. SB solutions obtained with $A$ and with $\tilde{A}_k$ at convergence for \Cref{Fig:blurred hst} using $\lambda=0.0500$ and $\beta=0.0008$. Here the PSF has size $11 \times 11$ pixels. }\label{Fig:hst} 
 \end{figure}

 \begin{exmp}\label{ex:hst}
 We consider the \texttt{hst} image of size $200 \times 200$ from \cite{gazzola2019ir}, shown in \Cref{Fig:True hst}. We blur this image with a mild non separable speckle PSF given in \Cref{Fig:PSF hst}. The provided PSF is a spatially invariant blur. Then we use $3 \%$ white Gaussian noise, corresponding to signal to noise ratio (SNR) of approximately $30$, to display the noise contaminated \texttt{hst} image in \Cref{Fig:blurred hst}. The parameters that define the blurring operator matrix and the number of terms $k$ used for the $\tilde{A}_k$ with the oversampling parameters $p$, as determined by the tolerance for $\nu_j $, are given in \Cref{Tab:Blur hst}.
\end{exmp}
\begin{center}
\begin{table}[ht!]
\centering
\caption{\label{Tab:Blur hst}The parameters that define the blurring matrix in \Cref{ex:hst} and the number of terms $k$ used for the $\tilde{A}_k$, with the oversampling parameter $p$, as determined by the tolerance for $\nu_j$.}
\begin{tabular}[t]{ c c c c c c c c c c }
\toprule
Blurring Mat &Size &Rank & Condition number& Kronecker rank & $k$ & $p$\\ \midrule
$A$ &$40000 \times 40000$& $37915$ & $4.1 \times 10^{19}$& $199$ & $16$& $4$ \\
\bottomrule
\end{tabular}
\end{table}  
\end{center}
The size of the problem is large, so we focus on obtaining the results using  SP. The errors in $\bar{\mathcal{R}}_k(A)$ and $\tilde{A}_k$ obtained for the matrix $A$ using EGKB(SP) and RSVD(SP) in \Cref{ex:hst} are presented in \Cref{Tab:Blur ex3}. 
 \begin{center}
\begin{table}[ht!]
\centering
\caption{\label{Tab:Blur ex3}The relative errors for
$\bar{\mathcal{R}}_k(A)$ as an approximation to $\mathcal{{R}}(A)$,  and for 
$\tilde{A}_k$ as an approximation to $A$ in \Cref{ex:hst}. The presented results for $\tilde{A}_k$ with EGKB(SP) are obtained with full reorthogonalization.} 
\begin{tabular}[t]{c c c c }
\toprule \\[-3ex]
Method &$\frac{\|\mathcal{R}(A)-{\bar{\mathcal{R}}}_k(A)\|_F}{\|\mathcal{R}(A)\|_F}$&$\frac{\|A-\tilde{A}_{{k}}\|_F}{\|A\|_F}$\\ 
\midrule       
EGKB(SP)&$0.0034$ &$0.0039$ \\
RSVD(SP)&$0.0034$&$0.0039$ \\
\bottomrule
\end{tabular}
\end{table}  
\end{center}

The parameters in the SB are chosen based on the noise level in this problem. Notably, the impact on the image in this problem caused by the PSF is comparable to that in \Cref{ex:Pattern2}, but here, we impose smaller white Gaussian noise on the image. Thus, we set $\gamma=30$ to minimize the effect of the shrinkage formula. Since the used noise level here is smaller than that used in \Cref{ex:Pattern2},  we use a smaller $\lambda=0.0500$, to yield $\beta=0.0008$. We note that the tolerance $\tau_{SB}=10^{-2}$ is suitable for demonstrating the efficiency of the SB and leads to faster convergence. The reconstructed \texttt{hst} images using the $\tilde{A}_k$ 
are illustrated in \Cref{Fig:hst}, and the quantitative results are summarized in \Cref{Tab:parameters hst}.

\begin{table}[ht!]
 \footnotesize
  \caption{\label{Tab:parameters hst}The corresponding \text{RE} and \text{ISNR}, \cref{eq:RE,eq:ISNR}, respectively,  by the SB with $\tilde{A}_k$ obtained using EGKB(SP), and RSVD(SP)  applied to the problem with $\text{SNR}\approx 30$ in \Cref{ex:hst}. The matrix $A$ is of size $40000\times 40000$. All results are obtained with the SB parameters   $\lambda=0.0500$ and $\beta=0.0008$.}
 \begin{center}
 \begin{tabular}[t]{c c c c c c c c c c c } 
 \toprule 
 &&&&&\multicolumn{5}{c}{CPU Time}\\ \cline{6-10}\\[-2ex]
Method&$\ell_{\text{end}}$&$\iota_{total}$&$\text{RE}(\bfx)$&$ \text{ISNR}(\bfx)$ 
 &$\mathcal{{R}}(A)$&$\bar{\mathcal{R}}_k(A)$&Terms($\tilde{A}_k$)& SB&Total\\ \midrule
ANI EGKB(SP) &$6$&$209$&$0.17$&$4.26$&$6.85$&$4.37$&$0.0029$&$4.55$&$15.77$\\  
 ANI RSVD(SP)&$6$&$209$&$0.17$&$4.26$&$6.85$&$3.94$&$0.0041$&$5.08$&$15.78$\\ 
 ISO EGKB(SP)&$2$&$72$ &$0.15$&$5.36$&$6.85$&$4.37$&$0.0029$&$1.71$&$12.93$\\ 
 ISO RSVD(SP)&$2$&$72$ &$0.15$&$5.36$&$6.85$&$3.94$&$0.0041$&$1.80$&$12.59$\\ 
\bottomrule 
 \end{tabular}       
   \end{center} 
 \end{table} 

Our presented results show that our approach to empirically estimate suitable choices for $\lambda$ and $\gamma$  to use in the SB algorithm, dependent on the noise level and the type of blur, is appropriate. For this problem, it is notable that  ISO SB yields a smaller relative error and requires approximately one third of the runs of $\ell_{\text{end}}$ and $\iota_{total}$ for the ANI algorithm, as indicated in \Cref{Tab:parameters hst}. As the problem size increases, the cost of the SB algorithm itself becomes negligible in comparison to forming the $\tilde{A}_k$, as is predicted in \Cref{Tab:costs}, when the number of iterations $\ell_{\text{end}}$ and $k$ are both small in contrast to the size of the problem. The RSVD(SP) is faster than EGKB(SP) in forming TSVD approximation to $\mathcal{R}(A)$, as reported in column $7$ of \Cref{Tab:parameters hst}. This result may contradict the costs, given in \Cref{Tab:costs}, and the results of the earlier experiments, presented in \Cref{ex:Pattern2} and \Cref{ex:satellite}. However, we expect that the use of full reorthogonalization impacts the performance of EGKB(SP) in this particular problem. The experiment illustrates that combining the approximation $\tilde{A}_k$ of $A$ with the SB algorithm can generate an acceptable restoration of the blurred image. 

In summary, \Cref{ex:Pattern2}, \Cref{ex:satellite}, and \Cref{ex:hst} demonstrate the performance of implementing the SB algorithm using the $\tilde{A}_k$ obtained using SP computations, as discussed in \cref{sec:Approx A} and \cref{subsec:SB}, respectively, for problems generated with various non-separable speckle PSFs and data contaminated with different noise levels.

\section{Conclusions}\label{sec:Conclusion}
In this paper we presented the use of the RSVD and EGKB algorithms to construct the Kronecker product sum approximation $\tilde{A}_k$ to $A$. The $\tilde{A}_k$ are demonstrated to be efficient in the context of image deblurring using anisotropic or isotropic total variation regularization implemented with the SB algorithm.
Results of numerical experiments on image deblurring problems illustrate the effectiveness of using the EGKB and RSVD to find $\tilde{A}_k$. Moreover, the estimate for $k$ to be used for the number of terms in the $\tilde{A}_k$ is done automatically within the EGKB algorithm using a new stopping test based on known theoretical results. Furthermore, our experiments demonstrate that we can implement the KP sum approximation portion of the algorithm using either RSVD or EGKB in single precision, allowing the application of the algorithm for larger matrices, while maintaining the quality of the results. Given the success of this approach, an important next step is determining a memory efficient reordering of the elements of $A$ when finding $\mathcal{R}(A)$ which may be a limiting factor for much larger implementations. 

The presented test using $\nu_j$ for stopping the EGKB algorithm, as described in \cref{subsec:testEGKB}, is a topic for further evaluation when finding a good low rank approximate TSVD for $A$ in other contexts. Further examination of the use of single precision in place of double precision for finding an approximate SVD in practical applications is also anticipated. 

\section*{Acknowledgments}
Funding: This work was partially supported by the National Science Foundation (NSF) under grant  DMS-2152704 for  Renaut.  Any opinions, findings, conclusions, or recommendations expressed in this material are those of the authors and do not necessarily reflect the views of the National Science Foundation.


\end{document}